\documentclass{article}
\usepackage[left=1.5in,right=1.5in,top=1.25in,bottom=1.25in]{geometry}

\usepackage{amssymb,amsmath,graphicx,comment,hyperref}

\newcommand {\grad}{\nabla}

\newcommand{\matlab}{{\sc Matlab}}

\newtheorem{theorem}{Theorem}

\newcommand{\beq}{\begin{equation}}
\newcommand{\eeq}{\end{equation}}
\newcommand{\beqf}{\begin{flalign}}
\newcommand{\eeqf}{\end{flalign}}

\usepackage{subfig}
\usepackage{float}
\usepackage{amsmath}
\usepackage[font=small,labelfont=bf]{caption}
\usepackage{caption}
\usepackage{soul}
\usepackage{cancel}
\usepackage{color}

\newcommand{\tr}[1]{\mathop{\mbox{Tr}}\left({#1}\right)}
\newcommand{\Diag}[1]{{\rm Diag}\left({#1}\right)}
\newcommand{\diag}[1]{{\rm diag}\left({#1}\right)}
\newcommand{\R}{{\mathbb{R}}}

\newcommand{\A}{\mathcal{A}}
\newcommand{\X}{\mathcal{X}}

\newcommand{\green}[1]{\textcolor{green}{#1}}

\newcommand{\cyan}[1]{USE MO MACRO}

\newcommand{\Hzero}{H_{k+1}^{(0)}}

\title{Behavior of Limited Memory BFGS when Applied to Nonsmooth Functions and their Nesterov Smoothings}
\author{Azam Asl\thanks{Booth School of Business, University of Chicago}
	\and 
	Michael L.~Overton\thanks{Courant Institute of Mathematical Sciences, New York University}
}	
\begin{document}
\maketitle

\begin{abstract}
The motivation to study the behavior of limited-memory BFGS (L-BFGS) on nonsmooth optimization problems 
is based on two empirical observations: the widespread success of L-BFGS  in solving large-scale smooth optimization 
problems, and the remarkable effectiveness of the full BFGS method in solving small to medium-sized nonsmooth optimization problems,
based on using a gradient, not a subgradient, oracle paradigm. We first summarize
our theoretical results on the behavior of the scaled L-BFGS method with one update applied to a simple convex
nonsmooth function that is unbounded below, stating conditions under which the method converges to a non-optimal point
regardless of the starting point. We then turn to empirically investigating whether the same phenomenon holds more generally,
focusing on a difficult problem of Nesterov, as well as eigenvalue optimization problems arising in semidefinite programming applications. 
We find that when applied to a nonsmooth function directly,  L-BFGS, especially its scaled variant, often breaks down with a poor
approximation to an optimal solution, in sharp contrast to full BFGS. Unscaled L-BFGS is less prone to breakdown but conducts far more function evaluations per iteration than scaled L-BFGS does, and thus it is slow. Nonetheless, it is often the case that both variants 
obtain better results than the provably convergent, but slow, subgradient method.
On the other hand, when applied to Nesterov's smooth approximation of a nonsmooth function, 
scaled L-BFGS is generally much more efficient than unscaled L-BFGS,
often obtaining good results even when the problem is quite ill-conditioned. 
Summarizing, we find that although L-BFGS is often a reliable method for minimizing ill-conditioned smooth problems,
when the condition number is so large that the function is effectively nonsmooth, L-BFGS frequently fails. 
This behavior is in sharp contrast to the behavior of full BFGS, which is consistently reliable for nonsmooth optimization problems. 
We arrive at the conclusion that, for large-scale nonsmooth optimization problems for which full BFGS and other methods for nonsmooth
optimization are not practical, it is often better to apply L-BFGS to a smooth approximation of a nonsmooth problem than to apply it 
directly to the nonsmooth problem.
\end{abstract}

\section{Introduction}  \label{sec:intro}

We consider the unconstrained optimization problem
	\[ \min_{x \in \R^n} ~~ f(x), \]
where the function $f$ is convex but nonsmooth. By this we mean that it is not differentiable everywhere,
and, typically, is not differentiable at minimizers.

Classical approaches to optimization of convex nonsmooth functions generally require the method to have access to an oracle that,
given $x\in\R^n$, returns the function value $f(x)$ and a subgradient  $g\in\partial f(x)$. The oldest such method, the subgradient
method of Shor, which dates to the 1970s, uses the iteration
\beq    \label{subgradmeth}
       x_{k+1} = x_k - t_k g_k, \text{   for some   } g_k\in\partial f(x_k),
\eeq
where $\{t_k\}$ is a pre-determined sequence of positive stepsizes. 
One well-known result states that, assuming $f$ is convex and bounded below, and provided the steplengths $\{t_k\}$
are square-summable (that is, $\sum_{k=0}^\infty t_k^2 < \infty$, and hence the steps are ``not too long"),
but not summable (that is, $\sum_{k=0}^\infty t_k = \infty$, and hence the steps are ``not too short"),
then convergence of $f(x_k)$ to the minimal value of $f$ must take place \cite{NB01}.
However, despite the strength of this theoretical result, it is well known that convergence to the optimal value is often very slow.
This observation led to the development of other algorithms, particularly the
bundle methods pioneered by Lemar\'{e}chal \cite{LEM75} for nonsmooth
convex functions in the mid-1970s, as well as the bundle methods of Kiwiel \cite{KIW85}
 for nonsmooth, nonconvex problems in the 1980s. 
As suggested by the name, at
each iteration, bundle methods use subgradient information obtained at former iterates as well as at the current iterate $x_k$ to
generate the next iterate $x_{k+1}$. Convergence results
are available for these methods, but the computational cost per iteration is significant for large $n$ as computing
$x_{k+1}$ from $x_k$ usually requires the solution of a quadratic program in $n$ variables \cite[p.\ 306 and 313]{BKM14}.
For more methods for nonsmooth optimization using the subgradient oracle, see \cite{BKM14,BGKMT20}.

In this paper, we do not use the subgradient oracle paradigm. One reason for this is that, in practice, in the
presence of rounding errors, it is often difficult, if not impossible, to determine whether the function $f$ is differentiable at a given point $x$
and hence to return a vector $g$ which can be guaranteed to be a subgradient.
A second is that since convex functions (and more generally, locally Lipschitz functions) are differentiable almost everywhere,
there is no reason, at least in the absence of exact line searches, to suppose that a method will ever generate a point $x_k$
where $f$ is not differentiable. We therefore use the simpler paradigm that, given $x$, an oracle returns $f(x)$ and $g = \grad f(x)$,
if $f$ is differentiable at $x$.  Clearly, the subgradient and gradient oracles coincide when $f$ is indeed differentiable
at $x$. What if this is not the case? In theory, we need to assume that the oracle informs the method that the function is not
differentiable and therefore a gradient cannot be provided. But in practice, the gradient oracle is
simply implemented as if $f$ \emph{is} differentiable at $x$, breaking ties arbitrarily if necessary.
 For example, if $f(x)$ is defined as $\max(f_1(x),f_2(x))$, where $f_1$ and
$f_2$ are smooth functions, the oracle returns $\grad f_1(x)$ if $f_1(x)>f_2(x)$, $\grad f_2(x)$ if $f_1(x)<f_2(x)$, and either one if
$f_1(x)=f_2(x)$, a property that is difficult to determine anyway in the presence of rounding errors. A subgradient oracle might,
in principle, return any convex combination of $\grad f_1(x)$ and $\grad f_2(x)$, but there is no reason for it to return anything other
than $\grad f_1(x)$ or $\grad f_2(x)$, unless, for example, the intent is to return a ``steepest descent" subgradient.

Using the gradient oracle instead of the subgradient oracle, we may ask the question: what can be said if we consider the
ordinary gradient method,
\beq  \label{gradmeth}
       x_{k+1} = x_k - t_k g_k, \text{   where   } g_k = \grad f(x_k),
\eeq
where the steplength $t_k$ is obtained, not from a predetermined sequence, but from a backtracking or Armijo-Wolfe line
search? The answer is well known: it is often the case that $x_k$ converges to a point $\bar x$
where $f$ is not differentiable and not minimal.  See \cite{AO20g} for a historical discussion and for a
detailed analysis of an interesting special case. Note that in relevant
illustrative examples, it is typical that $f$ \emph{is} differentiable 
at all iterates $\{x_k\}$, and nondifferentiable only at the limit point $\bar x$. This demonstrates that the power of Shor's subgradient
method is not so much that it is defined even if $f$ is not differentiable at $x_k$, as that the acceptable predetermined sequence of
stepsizes $\{t_k\}$, unlike stepsizes generated by a line search, ensures convergence to a minimal value.

In the early 2000s, Burke, Lewis and Overton introduced the gradient sampling algorithm \cite{BLO}
for nonsmooth, nonconvex optimization. This method uses the gradient, not the subgradient, oracle paradigm,
and a convergence analysis for a rather general class of functions states that, with probability one, the method
generates iterates $\{x_k\}$ where $f$ is differentiable and that,
if $f$ is convex and bounded below, the function values $\{f(x_k)\}$ converge to the minimal value
(more generally, that all cluster points of $\{x_k\}$ are Clarke stationary).
However, the cost per iteration is prohibitive
when the number of variables is large. See the survey paper \cite{BCLOS} for more details,
as well as recent work \cite{CurLi20} introducing more efficient variants of the gradient sampling method. 

As discussed by Lewis and Overton \cite{LO13}, the ``full" BFGS quasi-Newton method is a very effective alternative
choice for nonsmooth optimization, and its $O(n^2)$ cost per iteration 
is generally much less than the cost of the bundle or gradient sampling methods, 
but its convergence results for nonsmooth functions are limited to very special cases. 
It also uses the gradient, not the subgradient, oracle paradigm. 

Since the limited memory variant of BFGS \cite{LN89} (L-BFGS) costs only $O(n)$ operations per iteration, 
like the gradient and subgradient methods, it is natural to ask whether L-BFGS could be an effective method for
nonsmooth optimization. This is the topic of this paper, which is 
organized as follows. In \S \ref{sec:theory}, we summarize our theoretical results.
Then in \S \ref{sec:practice}, we give extensive experimental results. 
Although the methods we discuss are applicable to 
nonconvex problems, we restrict our discussion to the convex case for simplicity, focusing on
a difficult nonsmooth problem of Nesterov as well as eigenvalue optimization problems, including those arising
from semidefinite programming formulations of the max cut and matrix completion problems.
We also consider Nesterov's smooth approximations to these nonsmooth problems.
We make some concluding remarks in \S \ref{sec:conclu}.

\section{Limited Memory BFGS for Nonsmooth Optimization in Theory} \label{sec:theory}
We begin this section by defining the concept of an Armijo-Wolfe line search.
Then we discuss the full BFGS method on which L-BFGS is based, along with its
properties, before discussing our results for the L-BFGS method.

\subsection{Armijo-Wolfe Line Search}
The BFGS and L-BFGS methods rely on an Armijo-Wolfe line search (also often known as a ``weak Wolfe" line search).
Let $c_1$ and $c_2$ be parameters (the Armijo parameter and the Wolfe parameter respectively)
satisfying $0<c_1<c_2< 1$.  Let $x_k\in\R^n$ be a given iterate where $f$
is differentiable and let $d_k\in\R^n$ be a descent direction for $f$ from $x_k$, that is,
with ${\nabla f(x_k)}^Td_k < 0 $. We say that a positive steplength $t_k$ satisfies the Armijo condition if
\beq
f(x_k+t_k d_k) \leq f(x_k) + c_1 t_k \grad f(x_k)^T d_k,  \label{armijo}
\eeq
and that $t_k$ satisfies the Wolfe condition if $f$ is differentiable at $x_k+t_k d_k$ and  
\beq 
\grad f(x_k+t_k d_k)^T d_k \geq c_2 \grad f(x_k)^T d_k.  \label{wolfe}
\eeq
The Armijo condition imposes a ``sufficient decrease" in the value of $f(x_k+t_kd_k)$ compared to $f(x_k)$, while the Wolfe
condition imposes a ``sufficient increase" in the directional derivative $\grad f(x_k + t_k d_k)^T d_k$ compared to the negative
value  $\grad f(x_k)^T d_k$.  If $f$ is bounded below on the ray $\{x_k + td_k: t>0\}$, then it is known that points $t_k$ satisfying
these conditions exist under various assumptions on $f$, such as $f$ being convex (but not necessarily smooth).
For further discussion, including specification of a bracketing line search algorithm based on bisection and doubling to
compute $t_k$, see \cite{LO13}.

While it is often recommended in the context
of general smooth nonlinear optimization to set the Armijo parameter to a rather small value such as $10^{-4}$ \cite[p.\ 62]{NW06}, we note that for satisfactory complexity results for the gradient method applied to smooth, strongly convex functions, the Armijo parameter must \emph{not} be too small \cite[p.466--468]{BV}. 
An analysis of the gradient method using an Armijo-Wolfe line search applied to a nonsmooth
function that we discuss later in this paper (see \eqref{fdef_unbd}) was given in \cite{AO20g}.
It was shown that, in this case, the success or failure of the gradient method depends critically
on the Armijo parameter being sufficiently small, and experiments applying the gradient method to 
another nonsmooth function \eqref{leshouches} confirm the importance of
this issue \cite[Fig.\ 6]{AO20g}. On the other hand, as we will see in
\S\ref{subsec:lbfgs}, the choice of Armijo parameter is less critical to the success or failure of L-BFGS when applied
to the same function \eqref{fdef_unbd}.

\subsection{Full BFGS}\label{subsec:bfgs}

The  Broyden-Fletcher-Goldfarb-Shanno (BFGS) method, proposed independently by these four authors in 1970,
is a quasi-Newton method that maintains an approximation  $H_k$ to the  inverse of the  Hessian $\nabla^2 f(x_k)$.
Let an initial iterate $x_0 \in \R^n$ and an initial positive definite matrix $H_0$ be given.
The BFGS method is defined by the following iteration. For $k=0,1,2,\ldots$, let
\begin{align}
        d_k &= -H_k \nabla f(x_k) \nonumber \\
        x_{k+1} & = x_k + t_k d_k \text{  where  }t_k > 0 \text{ satisfies }\eqref{armijo},\eqref{wolfe}  \nonumber \\
        s_k & = t_k d_k   \label{skdef} \\
       y_k &  = \nabla f(x_{k+1}) - \nabla f(x_k) \label{ykdef}\\  
       \rho_k & = \frac{1}{s_k^T y_k}  \nonumber \\ 
      H_{k+1} & =  \left( I - \rho_k s_{k}y_{k}^T  \right)  H_{k}  \left( I -  \rho_k y_{k}s_{k}^T  \right) + \rho_k  s_{k}s^T_{k}  \label{bfgs_update} 
\end{align}
Equation \eqref{bfgs_update} is called the BFGS update. 
It is easy to check that the Wolfe condition ensures that $s_k^T y_k > 0$, and it is well known that the positive definiteness
of $H_{k+1}$ follows as a consequence. It is also clear that $H_{k+1}$ can be computed in $O(n^2)$ operations.

The most important convergence property of BFGS is due to Powell \cite{Pow76b}. It states that if $f$ is twice continuously
differentiable and strongly convex on the level set $S = \{x: f(x) \leq f(x_0)\}$, then the sequence $\{x_k\}$
converges to the global minimizer. Furthermore, the convergence theory of Dennis and Mor\'e  states that
if we further assume that the Hessian of $f$, $\grad^2 f(x)$, is Lipschitz continuous at the global minimizer, 
then the rate of convergence is superlinear.  See \cite{NW06} for more details.

For many years, BFGS was not considered as a possible stand-alone method for nonsmooth optimization, although
Lemar\'echal briefly mentioned this possibility in a 1982 technical report. Instead, Lemar\'echal
and others focused on using the BFGS update to provide second-order information to bundle methods \cite[p.~313]{BKM14}.
In contrast, Lewis and Overton \cite{LO13} advocated using the pure BFGS method, with no bundle method component,
as a practical method for nonsmooth optimization, using the gradient, not the subgradient, oracle paradigm. 
The Wolfe condition \eqref{wolfe} ensures in theory
that $f$ is differentiable at all iterates, although in practice, as already noted, it is difficult to check whether or not $f$ is differentiable 
at a given point. If $f$ is not differentiable at $x$, then the gradient oracle may return very different gradients
at points that are close to $x$. The consequence is that, in the nonsmooth case, BFGS builds a very ill-conditioned inverse ``Hessian'' 
approximation, with some tiny eigenvalues corresponding to ``infinitely large'' curvature in the directions
defined by the associated eigenvectors. This phenomenon, illustrated in \cite[Fig.\ 4]{LO13}, is apparently what makes
BFGS so effective for nonsmooth optimization.
Remarkably, the condition number of the inverse Hessian approximation often reaches 
$10^{16}$ before the method breaks down numerically, usually because, as a consequence of rounding error, the line search	
is unable to return $t_k$  satisfying the Armijo condition even though, in principle, an Armijo-Wolfe step exists.
Convergence of $f(x_k)$ to the minimal value of $f$ often appears to be linear (not superlinear, as in the smooth case).

Although empirically BFGS works very well when $f$ is nonsmooth, convergence results are limited to a few
special cases. The following results hold for all $x_0$ as long as $f$ is differentiable at $x_0$ and for all positive definite $H_0$.
We use $x^{(i)}$ to denote the $i$th entry of the vector $x$.

\begin{itemize} 
\item $n=1$  with $f(x)=|x|$:  the sequence generated by BFGS converges linearly to zero and is
is related to a certain binary expansion of the starting point \cite{LO13}.

 \item $f(x) = |x^{(1)}| + \sum_{i=2}^n x^{(i)}$: eventually a direction is identified on which $f$ is unbounded below
 \cite{XW17}; see also \cite{LZ15}.

 \item $f(x) = \|x\|_2$: iterates converge to the origin \cite{GL18}. This is a special case of a more general result 
 whose proof is based on Powell's theorem mentioned above.

\end{itemize}
As far as we know, even the case $f(x) = |x^{(1)}| + (x^{(2)})^2$ remains open!

BFGS has been used successfully in many practical applications of nonsmooth optimization including the
design of fixed-order controllers for linear dynamical systems with input and output,
shape optimization for spectral functions of Dirichlet-Laplacian operators and
condition metric optimization. For more details, see \cite[p.\ 159--160]{LO13}.
Software is available in the {\sc hanso} package.\footnote{http://www.cs.nyu/overton/software/hanso/}

Finally, BFGS has also proved useful for optimization problems with nonsmooth constraints.
Consider problems of the form 
\begin{align*}
    \min ~ &  f(x) \\ 
    \mathrm{~subject~to~} & c_i(x) \leq 0, ~~ i=1,\ldots,p  
\end{align*}
where $f$ and $c_1,\ldots,c_p$ may not be differentiable at local minimizers.
A successive quadratic programming (SQP) method based on BFGS was introduced by \cite{CMO17} 
to solve problems of this form, and applied to problems in static output feedback control design that involve 
spectral radius and pseudospectral radius constraints.
Although there are no theoretical results, it is typically much more efficient in practice than an SQP gradient sampling
method \cite{CurOve12} which does have convergence results. 
Software is available in the {\sc granso} package.\footnote{http://www.timmitchell.com/software/granso/}

\subsection{Limited Memory BFGS} \label{subsec:lbfgs}

Unlike the full BFGS method, L-BFGS does not store a matrix approximating the inverse of the Hessian $\grad^2f(x)$.
Instead, it uses an implicit approximation. Let $m$ be a small integer.
Instead of using \eqref{bfgs_update} to update a stored matrix $H_k$, 
the matrix $H_{k+1}$ is implicitly defined by applying $m$ BFGS updates successively,
starting with an initial matrix $\Hzero$ and using the updates defined by the pairs $(s_j,y_j), j=k-m+1,\ldots k$, defined in 
\eqref{skdef} and \eqref{ykdef}. This method is called L-BFGS-$m$. We consider two variants. In the \emph{scaled} variant, with
\beq
                \Hzero = \frac{s_k^T y_k}{y_k^T y_k} I,  \label{H0scaled}
\eeq
a scaling of the identity matrix often known as Barzilai-Borwein scaling \cite{BB88} is used
to initialize the implicit updating at iteration $k$.
In the \emph{unscaled} variant, $\Hzero$ is set to the identity matrix.
We also refer to these two different variants of L-BFGS-$m$ as \emph{with} and \emph{without scaling}, respectively.
Substantial numerical experience with applying L-BFGS-$m$ to minimize smooth functions
shows that the use of scaling is strongly preferred.  For more details on L-BFGS, including efficient implementation,
see \cite{LN89,NW06}.

Our theoretical analysis of the behavior of limited memory BFGS uses the scaled variant with $m=1$: $H_{k+1}$ is defined
by applying just one BFGS update to \eqref{H0scaled}. This method, L-BFGS-1, is sometimes called \emph{memoryless BFGS} 
\cite[p.\ 180]{NW06}.  We focus on the nonsmooth function
\beq
           f(x) = a|x^{(1)}| + \sum_{i=2}^n x^{(i)}   \label{fdef_unbd}
\eeq
with $a\geq \sqrt{n-1}$.   This function is obviously unbounded below, but it is bounded below
along any steepest descent direction $d=-\grad f(x)$.  One advantage of studying this function is its simplicity,
but another is that it is easy to determine whether a method succeeds or fails when it is applied to \eqref{fdef_unbd}:
a method is successful only if it generates a sequence of function values $f(x^{(k)}) \rightarrow -\infty $ or identifies
a direction $d_k$ on which $\{f(x_k + t d_k): t > 0\}$ is unbounded below. Otherwise, since Armijo-Wolfe steps
always exist along directions on which $f$ is bounded below, the sequence of function values $f(x^{(k)})$ is well
defined, bounded below, and must converge to a non-optimal finite value. The following theorem shows that
this last scenario occurs when the Armijo parameter is not sufficiently small compared to the parameter $a$ defining
the function definition in \eqref{fdef_unbd}. The proof is given in \cite[\S 3.2]{AO20b}.

\begin{theorem}
Suppose $f$ is defined by \eqref{fdef_unbd} with $a \ge 2\sqrt{n-1}$.
Set $x_0$ to \emph{any} point with $x_0^{(1)}\not = 0$.
If the  Armijo parameter $c_1$ is chosen so that
		$$ \frac{1-c_1}{c_1}(n-1) < a^2+a\sqrt{a^2 - 3(n-1)}$$
holds, then the scaled L-BFGS-1 method is well defined, in the sense that Armijo-Wolfe steps always exist,
but fails in the sense that $f(x_k)$ is bounded below as $k \to \infty$. 
\end{theorem}

Note that $a \ge 2\sqrt{n-1}$ implies
\[
a^2+a\sqrt{a^2 - 3(n-1)} \geq 4(n-1) + 2\sqrt{n-1}\sqrt{n-1} = 6 (n-1).
\]
So, if  $c_1 > 1/7$, the failure condition holds.  
It was proved in \cite{AO20g} that, when applied to  \eqref{fdef_unbd},
 the gradient method with an Armijo-Wolfe line search fails in the same sense if the stronger condition
$$ \frac{1-c_1}{c_1}(n-1)  < a^2$$
on the Armijo parameter holds. So, surprisingly, in some cases scaled L-BFGS-1 fails although the gradient method succeeds
in generating $f(x_k)\rightarrow -\infty$.

If the specific Armijo-Wolfe bracketing line search given in \cite{LO13} is used, we obtain a failure condition for scaled L-BFGS-1
that is \emph{independent} of the Armijo parameter. The proof of the next result is given in \cite[\S 3.3]{AO20b}.

\begin{theorem}
Suppose $f$ is defined by \eqref{fdef_unbd} with $a \ge 2\sqrt{n-1}$. 
Set $x_0$ to \emph{any} point with $x_0^{(1)}\not = 0$.
If scaled L-BFGS-1 is implemented using the Armijo-Wolfe bracketing line search of \cite{LO13},
the method is well defined in the sense that the line search always returns an Armijo-Wolfe steplength $t_k$,
but fails in the sense that $f(x_k)$ is bounded below as $k \to \infty$. 
\end{theorem}

In fact, in the experiments reported in \cite{AO20b} we observe that $ a \ge \sqrt{3(n-1)} $ suffices for the method to fail.
Furthermore, experiments indicate that, for any number of updates $m$, there is a threshold for the parameter $a$
beyond which scaled L-BFGS-$m$ always fails in the same sense: $f(x_k)$ is bounded below as $k \to \infty$. This threshold value for $a$
increases with $m$.

\section{Limited Memory BFGS for Nonsmooth Optimization in Practice}\label{sec:practice}

In this section we carry out extensive experiments applying L-BFGS-$m$, with and without scaling, using the Armijo-Wolfe
bracketing line search of \cite{LO13}, to several challenging nonsmooth examples.  All of the functions we consider from now on are bounded below.
In most cases, in assessing whether a method succeeds or fails, we compare the final computed result with the known optimal
value. We also compare the results with those obtained by the full BFGS method using the same line
search. We used the {\sc hanso} package which implements both full and limited-memory BFGS, with the default  settings for the Armijo and Wolfe parameters: $c_1=10^{-4}$ and $c_2=0.5$. In many cases, we also compared
the results to those obtained by the subgradient method with predetermined stepsizes $t_k=1/k$. 
All methods, including the subgradient method, are implemented using
the gradient oracle paradigm discussed in \S \ref{sec:intro}: no attempt is made to determine whether the objective function is
differentiable at a given iterate $x_k$. 
Instead of using a stopping criterion such as discussed in \cite[p.\ 159]{LO13}, we run each method until it reaches a maximum number of function evaluations or iterations, or the line search fails to find an acceptable step (either because a limit on the number of bisections in the
line search is exceeded or because of rounding errors).  If the method terminates because of failure in the line search when
$f(x_k)$ is not a good approximation to the optimal value, we say that the method ``breaks down''.  
We find, in general, that this behavior is typical for L-BFGS, but not for full BFGS. For this reason we also 
compare these methods on \emph{smoothed} versions of the objective functions.

We begin with a difficult nonsmooth problem devised by Nesterov, and then we go on to consider eigenvalue optimization
and semidefinite programming problems. 

\subsection{Nesterov's Les Houches Problem}\label{subsec:ylh}

We now consider a nonsmooth function introduced by Nesterov.
The function is defined by
\beq\label{leshouches}
f(x) = \max\{|x^{(1)}|, ~|x^{(i)}-2x^{(i-1)}|,~~ i=2,...,n\}.
\eeq
Let $\hat x^{(1)}=1, \hat x^{(i)} = 2\hat x^{(i-1)} + 1, i=2,...,n$.
Then $f(\hat x) = 1 = f({\bf 1})$, where ${\bf 1}=[1,\ldots,1]^T$, although $\|\hat x\|_\infty \approx 2^n$ and $\|{\bf 1}\|_\infty=1$,
so the level sets of $f$ are very distorted.
The minimizer is ${\bf 0}=[0,\ldots,0]^T$ with $f({\bf 0})=0$. We call this Nesterov's ``Les Houches" function \cite{YN16}.

\begin{figure} 
	\includegraphics[width=5.5in,height=1.5in]{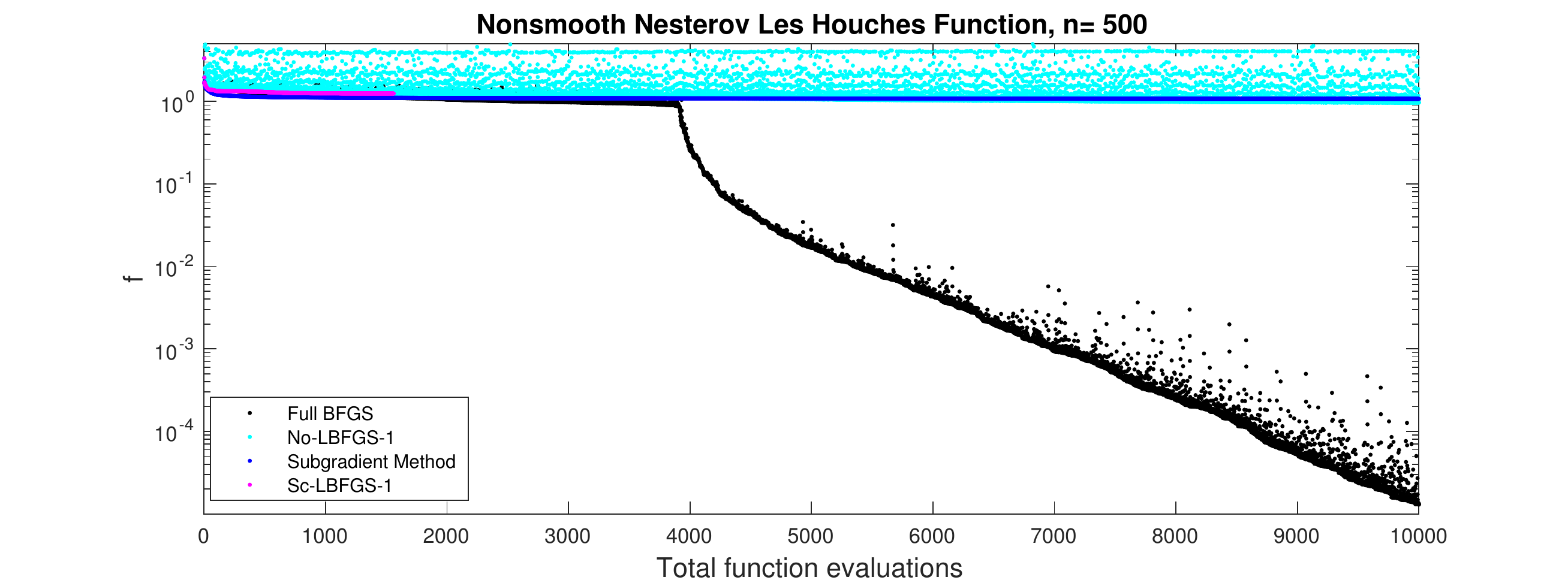} 
	\caption
	{Comparing full BFGS, L-BFGS-1 with and without scaling and the subgradient method on the nonsmooth Les Houches problem \eqref{leshouches} with $n=500$. Here and below, ``No-LBFGS" and ``Sc-LBFGS" refer to the
	methods without scaling and with scaling respectively. }
	\label{fig:ylh500}
\end{figure}

In Figure \ref{fig:ylh500}, we show results for L-BFGS-1 (memoryless BFGS), with and without scaling, as well as full BFGS and the subgradient method, for the Les Houches problem with $n=500$. We display all the function values that were computed, including those computed
in the bracketing line search. We put a limit of 10,000 function evaluations on each method. The starting point $x_0$ (used by all the methods) was  drawn randomly from the ball of radius $0.1$ centered at the vector of all ones, using the standard normal distribution.
\begin{figure} 
	\includegraphics[width=5.5in,height=1.5in]{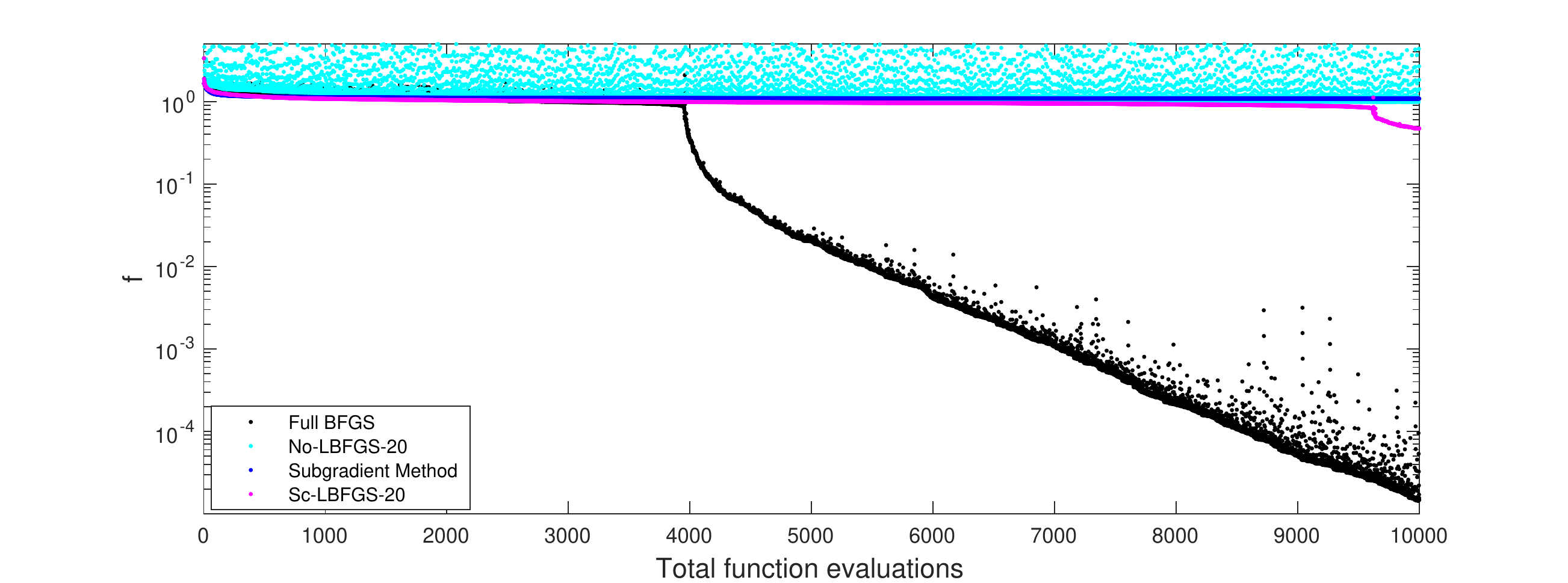}
	\caption
	{Comparing full BFGS, L-BFGS-20 with and without scaling and the subgradient method on the nonsmooth Les Houches problem \eqref{leshouches} with $n=500$. 
	}
	\label{fig:ylh500_20}
\end{figure}

We see from Figure \ref{fig:ylh500} that scaled L-BFGS-1 (magenta dots) breaks down, with failure in the line search, after fewer than 2000 function evaluations.  
In contrast, unscaled L-BFGS-1 (cyan) runs for the full 10,000 function evaluations. However, its scattered plot indicates that the method performs many function evaluations per iteration in the line search, indicating that, not surprisingly given its name, the search directions it generates are not well scaled. Despite this, the method obtains a somewhat lower answer than the subgradient method (dark blue). Full BFGS (black) performs much better than any of the other methods, reducing the function value to about $10^{-5}$ (recall that the optimal value is zero). It is interesting to note that its convergence rate picks up rapidly right after it has lowered the function value down to 1. We do not know the reason for this.

We now increase the number of updates $m$ from 1 to 20 and repeat this experiment: see Figure \ref{fig:ylh500_20}. Unlike scaled L-BFGS-1, scaled L-BFGS-20 does not quit early, and furthermore it also demonstrates a suddenly faster convergence rate toward the end of the experiment similar to that of full BFGS (although the final answer it obtains is not nearly as accurate as full BFGS). It obtains a  function value of size 0.47, whereas unscaled L-BFGS-20 gets a final answer of about 0.998 and the subgradient method obtains 1.083. In this experiment, increasing the number of updates from 1 to 20 enhanced the performance of scaled L-BFGS far more than it did for unscaled L-BFGS.
 
The conclusion from these experiments is that with a small $m$, unscaled L-BFGS-$m$ performs better and with a larger $m$ it is the scaled variant which performs better. However, neither method performs nearly as well as full BFGS.


\subsection{Smoothed Versions of Nesterov's Les Houches Problem}\label{subsec:sylh}
Since L-BFGS-$m$ performed poorly on \eqref{leshouches}, we consider instead applying it to a smoothed version.
Let
\[
A=
\begin{bmatrix}
1 & 0 & 0 & \hdots  & 0 \\
-2 & 1 & 0 & \hdots & 0 \\
\vdots & \vdots & \vdots & \vdots & \vdots \\
0 &  \hdots & 0  & -2 & 1
\end{bmatrix} .
\] 
Then, \eqref{leshouches} is equivalent to \beq\label{fgeq}
f(x)=g(Ax),
\eeq
where $g:\R^n \to \R$ is defined by
\beq\label{vec_max}
g(y) = \max\{|y^{(i)}| :i=1,2,...,n\}.
\eeq
Consider the Nesterov smoothing \cite{Nes05,VandenbergheCourse} of the vector-max  problem \eqref{vec_max}:
\beq \label{sdef}
g_{\mu}(y)= \mu\log\sum_{i=1}^n \big(e^{y^{(i)}/\mu} + e^{-y^{(i)}/\mu} \big) -\mu\log(2n).
\eeq
Without the constant term $-\mu\log(2n)$, this function is sometimes known as the softmax function \cite[p.\ 205]{Mohri12}. 
The unique minimizer of $g_{\mu}(y)$ is $y_\mu^*={\bf 0}$ with $g_{\mu}^*= 0$. Since $A$ is full-rank (although one of its singular values converges to zero as $n \to \infty$), via \eqref{fgeq} we know that the unique minimizer of $f_\mu(x)$, $x_\mu^*$, is also {\bf 0}, with the same optimal value $f_\mu^* = 0,$ and it can be verified that
\beq\label{hessbd}
\|\nabla^2 f_{\mu} ({\bf 0})\|_2  = \dfrac{1}{n\mu}\|A\|_2^2.
\eeq
We follow the standard approach \cite[Sec 9.1]{BV} to defining the condition number of the strongly
convex function $f_\mu$ as 
\beq\label{conddef}
\kappa(f_\mu) = \left (\max_{x\in S} \| \nabla^2 f_{\mu}(x)\|_2\right )
\left (\max_{x\in S} \|(\nabla^2 f_{\mu}(x))^{-1}\|_2\right)
\eeq
where $S = \{x: f(x) \leq f(x_0)\}$.  For small $\mu$,  the second factor is enormous as all eigenvalues of $\nabla^2 f_\mu(x_0)$ are tiny. Using \eqref{hessbd} as a lower bound for the first factor we conclude that, for small $\mu$,
\[
\kappa(f_\mu) \gg \frac{1}{\mu}.
\]

\begin{figure}
	\centering
	\begin{tabular}{cc}
	\includegraphics[scale=0.30]{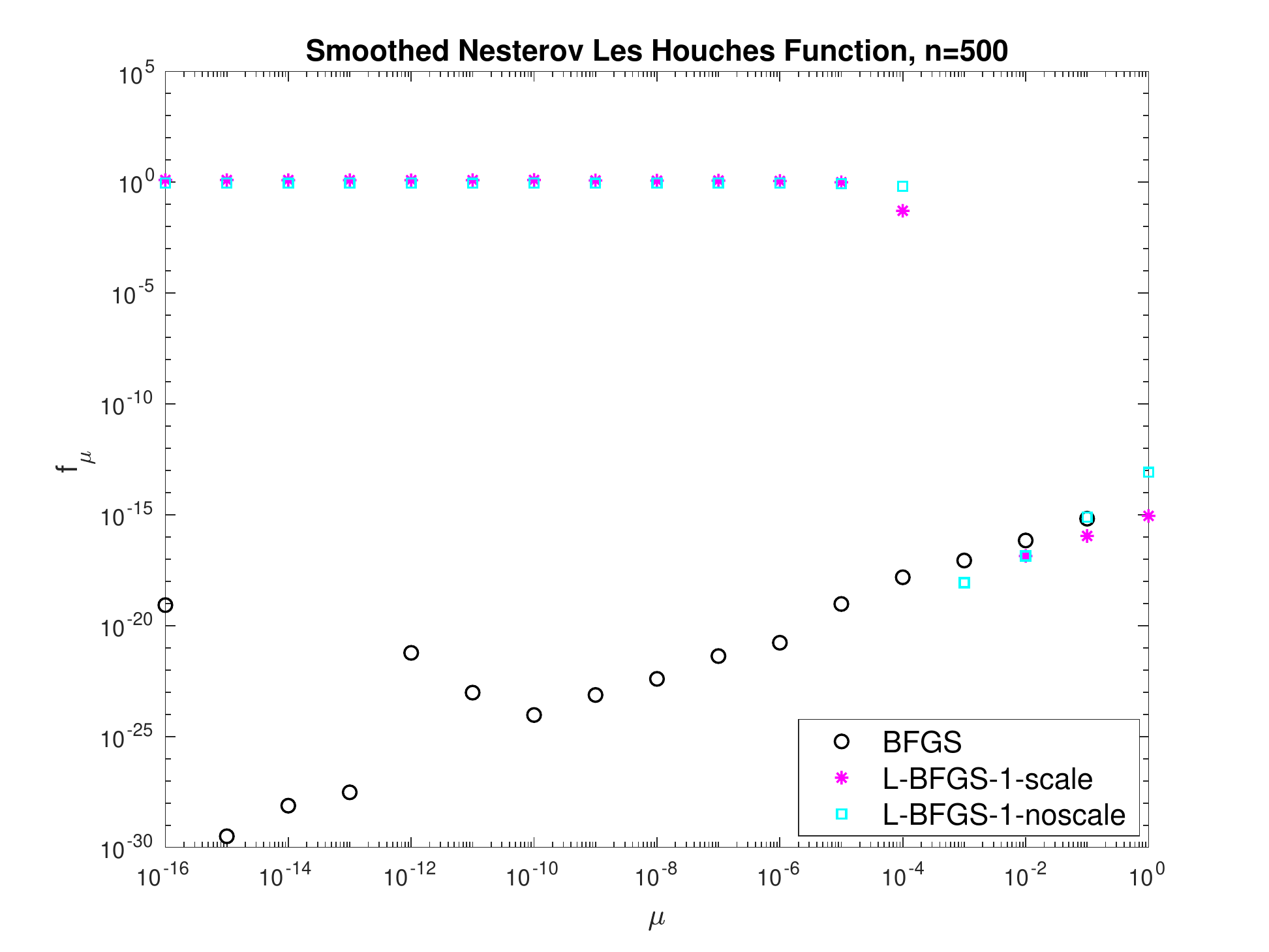} 
	&
	 \includegraphics[scale=0.30]{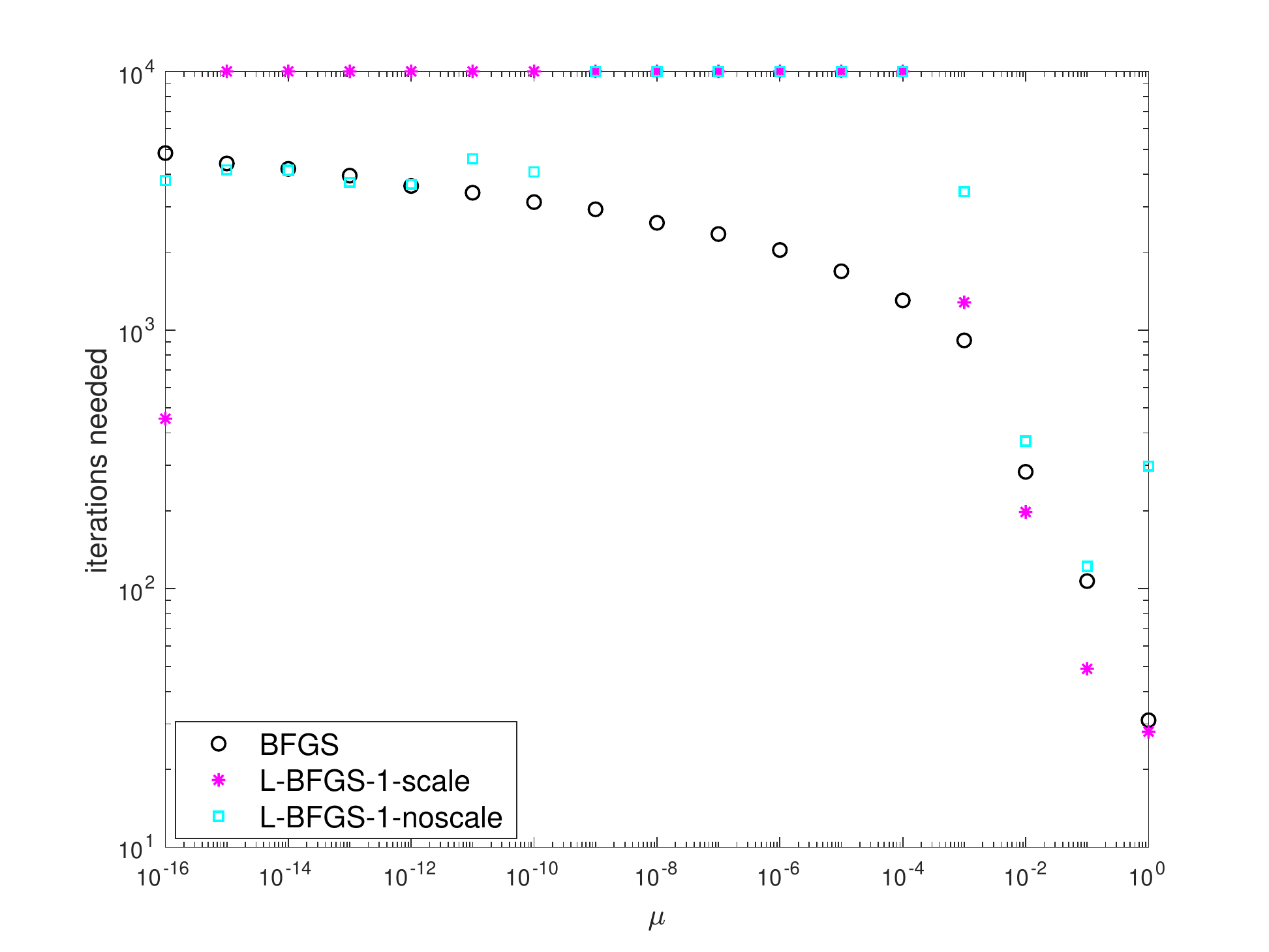}
	\end{tabular}
	\caption
	[ Smoothed Nesterov Les Houches  Function -L-BFGS-1]
	{Comparing BFGS and L-BFGS-1 with and without scaling on the smoothed Les Houches function \eqref {sdef} for $n=500$. The left panel shows the final function value and the right panel shows the iteration count, both as a function of the smoothing parameter $\mu$. The maximum number of iterations is set to $10^4$.}
	\label{fig:sylh1}
\end{figure}

We now report on experiments we conducted applying full BFGS and L-BFGS, with and without scaling, to the smooth function $f_\mu$, with $n=500$ as before. All of the methods start from the same initial point used earlier. 
The left panel of Figure \ref{fig:sylh1} shows the final function value computed by full BFGS and L-BFGS-1 with and without scaling as a function of the smoothing parameter $\mu$ using a log-log scale. Let us focus first on the results for full BFGS (black circles). 
BFGS always finds a solution with magnitude smaller than $10^{-15},$ even for a very small $\mu$, when the function is extremely ill conditioned. This is a remarkable property of full BFGS: its accuracy does not deteriorate significantly\footnote{Surprisingly, the accuracy increases somewhat as $\mu$ decreases, but this is at the level of rounding errors and could perhaps be explained by a rounding error analysis. Certainly the scatter at the bottom left corner of the left panel of Figure \ref{fig:sylh1} is a consequence of rounding error.}  as the condition number $\kappa(f_\mu)$ of the smoothed problem blows up with $\mu\rightarrow 0$. In fact, when $\mu$ is sufficiently small, say $\mu=10^{-16}$ (approximately the rounding unit in IEEE double precision used by {\sc{matlab}}), the smoothed problem is essentially equivalent to the original nonsmooth problem when rounding errors are taken into account, so the left panel shows the transition of the accuracy of full BFGS from smoothed variants of the problem to the limiting nonsmooth problem. The right panel shows the number of iterations that were required, again as a function of the smoothing parameter $\mu$ and again using a log-log scale. The maximum number of iterations (not function evaluations) was set to $10^{4}$ for each  $\mu$. Remarkably, we see that the number of iterations required for full BFGS to accurately minimize $f_\mu$ does not significantly increase as $\mu\rightarrow 0$, even though the condition number $\kappa(f_\mu)$ blows up as $\mu$ decreases to zero, and the number required for the effectively nonsmooth instance $\mu=10^{-16}$ is not much more than the number required for much better conditioned smoothed problems arising from moderate values of $\mu$.

The results for L-BFGS-1 are very different. Unscaled L-BFGS-1 (cyan squares) finds an accurate answer for $\mu \geq 10^{-3}$, but the number of iterations required increases rapidly as $\mu$ is decreased further so the iteration limit is reached for $\mu$ ranging from $10^{-4}$ to $10^{-9}$. However, starting with $\mu=10^{-10}$, unscaled L-BFGS-1 breaks down before reaching the maximum number of iterations. The behavior of scaled L-BFGS-1 (magenta asterisks) is similar except that it breaks down only for $\mu=10^{-16}$. Note that, since we are displaying the number of iterations, not the number of function evaluations, the performance of unscaled L-BFGS-1 looks better than it really is: the scaled version is computing substantially fewer function evaluations per line search.
\begin{figure}
    \centering
	\begin{tabular}{cc}
	\includegraphics[scale=0.3]{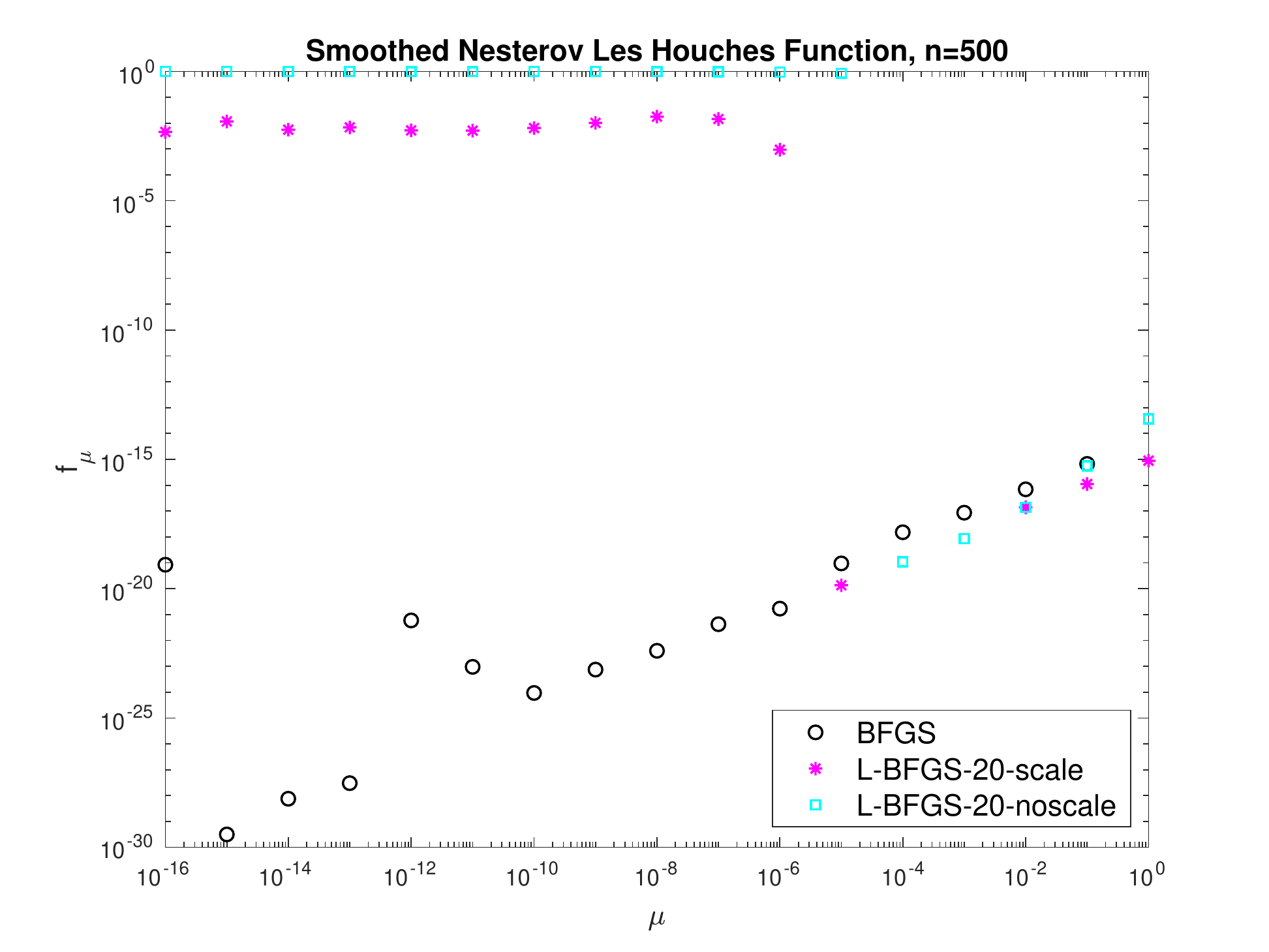} 
	&
	\includegraphics[scale=0.3]{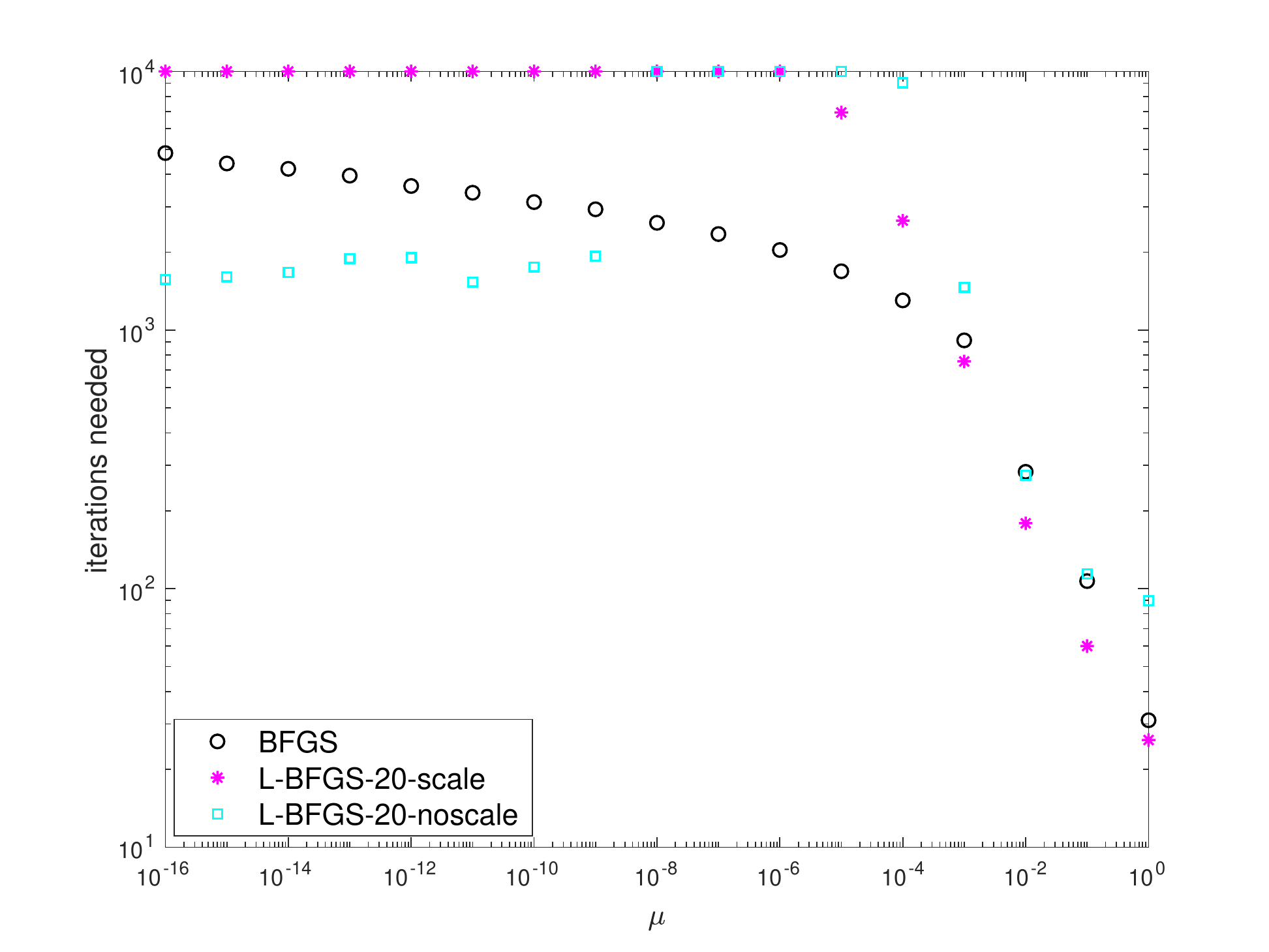} 
	\end{tabular}
	\caption
	[Smoothed Nesterov Les Houches  Function -L-BFGS-20]
	{Comparing BFGS and L-BFGS-20 with and without scaling on the smoothed Les Houches function \eqref {sdef} for $n=500$. The left panel shows the final function value and the right panel shows the iteration count, both as a function of the smoothing parameter $\mu$. The maximum number of iterations is set to $10^4$.}
	\label{fig:sylh20}
\end{figure}

When we increase the number of updates to $m=20$, scaled L-BFGS reacts much better than unscaled L-BFGS; see Figure \ref{fig:sylh20}. Unscaled L-BFGS-20 finds accurate answers for \mbox{$\mu\geq 10^{-4}$} before hitting the iteration limit, while the scaled version does so for $\mu\geq 10^{-5}$. Furthermore, when the maximum iteration limit is reached, scaled L-BFGS-20 achieves an answer of magnitude $\approx$ $10^{-2}$, whereas unscaled L-BFGS-20 is still giving an answer of magnitude $\approx$ $10^{0}$, similar to unscaled L-BFGS-1. However, overall, both are still doing poorly compared to full BFGS.

Our conclusions from this subsection are consistent with the generally accepted wisdom concerning L-BFGS. For smooth problems, even very ill-cond\-ition\-ed ones, it is best to use the scaled version of L-BFGS, and choosing the number of updates $m$ to be larger rather than smaller gives better performance, although, in contrast to full BFGS, the number of iterations required increases significantly with the conditioning of the problem. Again in contrast with full BFGS, when the ill conditioning increases to the nonsmooth limit implicit in consideration of rounding errors, scaled L-BFGS generally fails to converge to an optimal solution. However, for the smooth but ill-conditioned problems considered in this subsection, 
unscaled L-BFGS offers no advantage compared to scaled L-BFGS. The most important conclusion is that, while applying full BFGS directly to nonsmooth problems works remarkably well, this is not the case for L-BFGS; at least for the Les Houches problem, it is far preferable to apply scaled L-BFGS to a smooth approximation of the nonsmooth problem. 

\subsection{Max Eigenvalue Problem}
\label{subsec:max_eig}
Let $S^N$ denote the space of $N\times N$ real symmetric matrices, and let $\A:S^N\to \R^n$ denote a linear operator acting on $X$ as follows:
\beq \label{Alinopdef}
\A X = \left[\begin{array}{c} \langle A_1, X	\rangle \\
	\vdots \\
	\langle A_n, X	\rangle\end{array}\right],
\eeq
with $A_i \in S^N$ for $ i=1., \hdots, n$. Its 
adjoint operator; $\A^T: \R^n \to  S^N$, is defined by
\beq \label{Alinopdef2}
\A^Ty  = \sum_{i=1}^{n} y^{(i)}A_i.
\eeq
The Max Eigenvalue problem is to minimize the function
\beq \label{max_eig}
f(y) = \lambda_{\rm max}(C-\A^Ty),
\eeq
where $C\in S^N$ and $\lambda_{\rm max}: S^N \to \R$  denotes largest eigenvalue of its argument. It is well known
that $\lambda_{\rm max}$ is a convex function on $S^N$. Early papers on eigenvalue optimization include \cite{MO88}.

Assuming the maximum eigenvalue of $C-\A^Ty$ is simple, the gradient of $f$ is
\[
\nabla f(y) = -\A (qq^T) = -[q^T A_1 q , \cdots , q^T A_n q]^T\green{.}    
\]
Following the gradient oracle paradigm discussed in \S \ref{sec:intro}, we make no attempt to estimate whether or not the maximum eigenvalue
is simple at a given iterate.
However, at optimal solutions, we generally expect that $C-\A^Ty$ has a multiple largest eigenvalue and hence $f$ is not differentiable.  As is well known, eigenvalue optimization problems are instances of semidefinite programs, and hence small problems can be solved using CVX \cite{cvx}.

\begin{figure}
	\includegraphics[width=5.5in,height=1.5in]{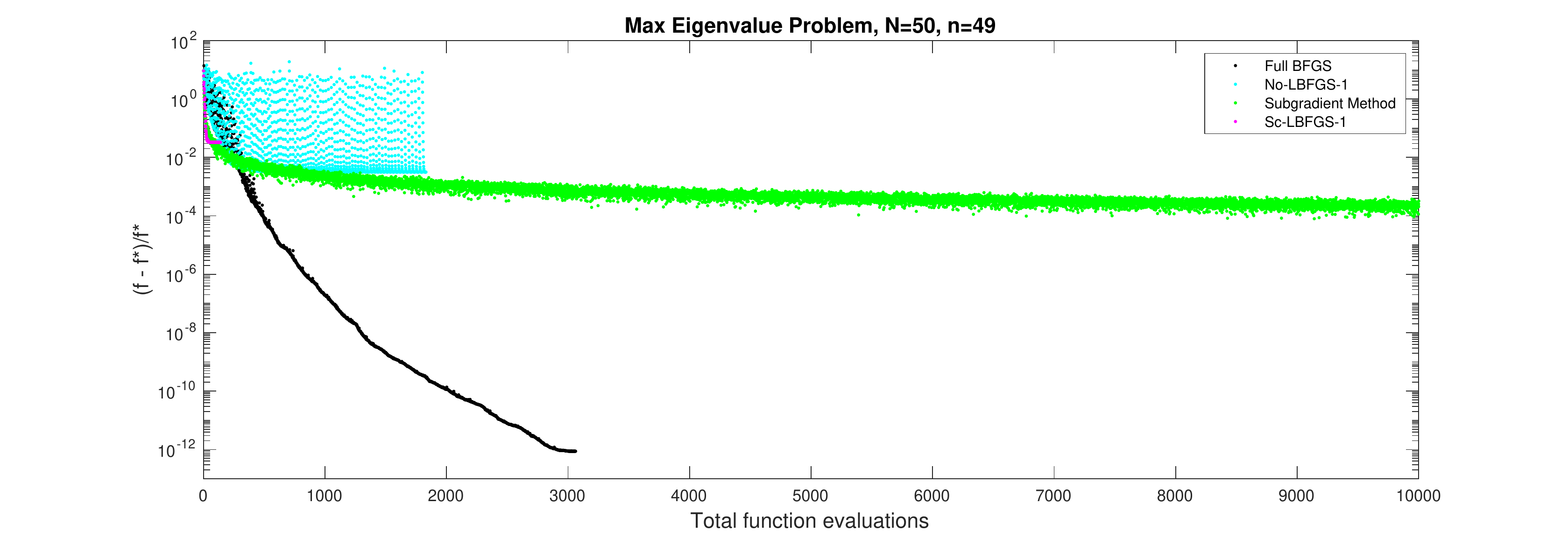}
	\caption
	[Random Max Eigenvalue Problem -L-BFGS-1]
	{Comparing BFGS, L-BFGS-1 with and without scaling and the subgradient method on a randomly generated Max Eigenvalue problem \eqref{max_eig} with $N=50$ and  $n=49$.}
	\label{fig:me1}
\end{figure}

Using the standard normal distribution, we generated a random instance of this problem, defining $\A$ and $C$  with  $N=50$ and $n=49$. Figure \ref{fig:me1} 
shows the performance of full BFGS, L-BFGS-1 with and without scaling, and the subgradient method (as before with $t_k=1/k$) for minimizing \eqref{max_eig}. As earlier, we display all the function values that were computed, including those computed
in the bracketing line search.
All methods were terminated after $10^4$ function evaluations. Each function evaluation $f(y)$ makes a call to the  \matlab\  function {\tt eig}  to compute all the eigenvalues of $C-\A^Ty$.
The vertical axis shows the relative error $|(f-f^*)/f^*|$, where we used the SDPT3 solver in CVX to obtain the optimal solution $f^*$ with accuracy $10^{-14}$. As in the experiment on the nonsmooth Les Houches problem reported in Figure \ref{fig:ylh500}, scaled L-BFGS-1 (magenta dots) breaks down early. However, unlike in that experiment, here  unscaled L-BFGS-1 (cyan dots) also breaks down early, and as a result the subgradient method (green dots) obtains a better answer, though not nearly as good as full BFGS (black dots).

It's also of interest to examine the multiplicity of the eigenvalues of $C-\A^Ty$ at the optimal solution $y^*$ and its computed approximations. From SDPT3, we know that the
optimal multiplicity for this problem is 5. Table \ref{table:mee} shows the top 6 eigenvalues of the final answer found by each method. Besides SDPT3, only BFGS and the subgradient method are able to determine the correct optimal multiplicity. BFGS finds a solution with 11 correct digits, while the subgradient method obtains 3 correct digits. Both variants of L-BFGS-1 converge to answers with multiplicity 4. Although the multiplicity is wrong, unscaled L-BFGS-1 gets a better answer than its scaled counterpart, with 2 correct digits, but at the cost of requiring
many more function evaluations. 

\begin{table}
         \centering
	\footnotesize
	\begin{tabular}{|c|c|c|c|c|}
		\hline
		SDPT3 & BFGS  & Sc L-BFGS-1 & No L-BFGS-1 & subgradient \\
		\hline
		7.82702970305352&  7.82702970306035&  8.08455876518360&  7.85155000878711&  7.82953885641783\\
		7.82702970305349&  7.82702970306035&  8.08455876518359&  7.85155000044960&  7.82746561454846\\
		7.82702970305348&  7.82702970306034&  8.08197715145863&  7.85154996050940&  7.82673229360627\\
		7.82702970305346&  7.82702970306031&  8.05541534062475&  7.85141043900471&  7.82472408900949\\
		7.82702970305334&  7.82702970306017&  7.84362627205676&  7.69075549655739&  7.82286893229481\\
		7.70350538019538&  7.70350432059448&  7.56258926523925&  7.48734455288558&  7.70188538367848\\
		\hline
	\end{tabular}
	\normalsize
	\caption
	[Top Eigenvalues of Random Max Eigenvalue Problem -LBFGS-1]
	{Top 6 eigenvalues of $C-\A^T y$ for Max Eigenvalue problem \eqref{max_eig} where $y$ is computed by SDPT3, full BFGS, scaled/unscaled L-BFGS-1   and the subgradient method for a randomly generated problem with $N=50$ and  $n=49$. The optimal multiplicity is 5. }
	\label{table:mee}
\end{table}

\begin{figure}
	\includegraphics[width=5.5in,height=1.5in]{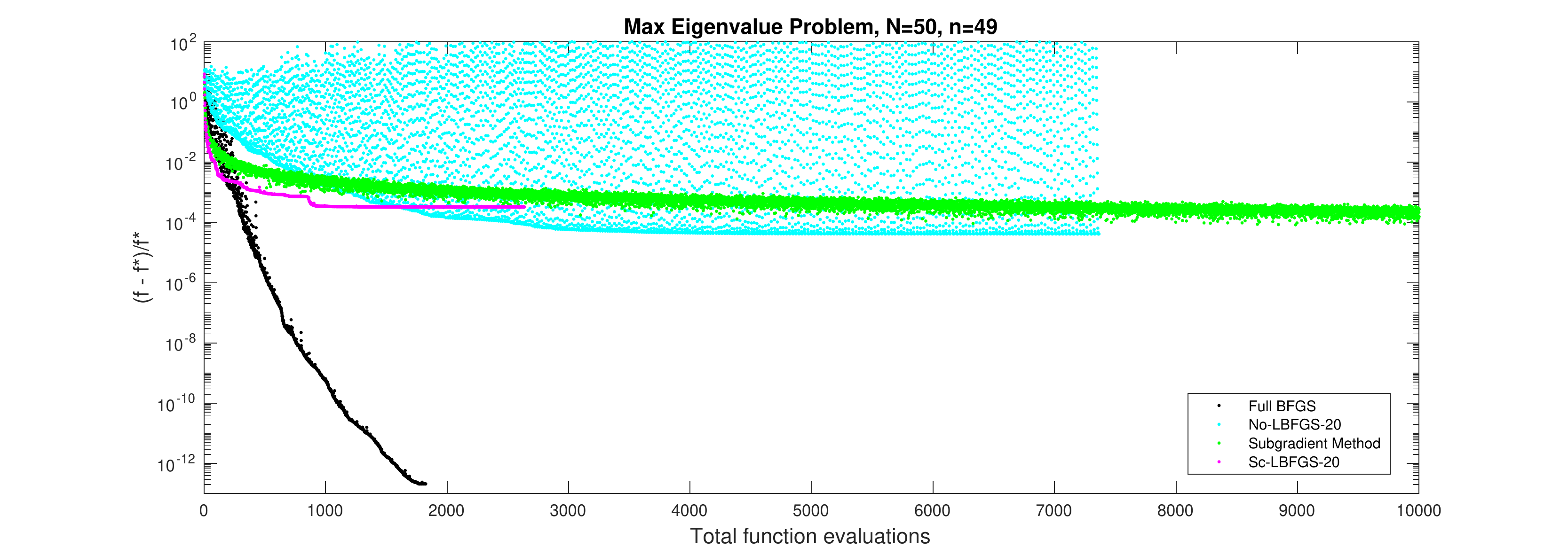} 
	\caption
	[Random Max Eigenvalue Problem -L-BFGS-20]
	{Comparing BFGS, L-BFGS-20 with and without scaling and subgradient method on a randomly generated Max Eigenvalue problem \eqref{max_eig} with $N=50$ and  $n=49$.}
	\label{fig:me20}
\end{figure}

Next, we  repeat this experiment on the same problem, increasing the number of L-BFGS updates from $m=1$ to $m=20$. See Figure  \ref{fig:me20} as well as Table \ref{table:mee20} which presents the top 6 eigenvalues for the final answer obtained by scaled and unscaled L-BFGS-20.  Both methods find the right multiplicity, with scaled L-BFGS-20 obtaining 3 correct digits and unscaled L-BFGS-20 obtaining 4 correct digits.
\begin{table}
	\centering
	\footnotesize
	\begin{tabular}{|c|c|c|c|c|}
                \hline
		Sc L-BFGS-20 & No L-BFGS-20 \\
		\hline
		7.82959659952176& 7.82735384547039\\
		7.82959659952176& 7.82735380191247\\
		7.82959659952174& 7.82735377961470\\
		7.82959659952166& 7.82735377236995\\
		7.82959659950328& 7.82735372682267\\
		7.62982269438813& 7.69181664817458\\
		\hline
	\end{tabular}
	\caption
	[Top Eigenvalues of Random Max Eigenvalue Problem -LBFGS-20]
	{Top 6 eigenvalues of $C-\A^T y$ for Max Eigenvalue problem \eqref{max_eig} where $y$ is computed by scaled and unscaled L-BFGS-20 for the same problem reported in Table \ref{table:mee}. The optimal multiplicity is 5.}
	\label{table:mee20}
\end{table}

In summary, we observe that for the Max Eigenvalue problem, unlike the  Les Houches problem, increasing  $m$ from 1 to 20 does not result in scaled L-BFGS doing better than unscaled L-BFGS.

\subsection{Smoothed Max Eigenvalue Problem}\label{subsec:smax_eig}

Consider now Nesterov smoothing of the Max Eigenvalue function \eqref{max_eig} \cite{AdA08} 
\beq\label{smax_eig}
f_{\mu}(y) = \mu\log \sum_{i=1}^N \exp( \lambda_i(C-\A^Ty)/\mu) - \mu\log N,
\eeq
where 
$\lambda_1(W)\geq \lambda_2(W) \geq \ldots \geq \lambda_N(W)$ denote the ordered eigenvalues of a symmetric matrix $W\in S^N$. Thus, $\lambda_1$ is equivalent to $\lambda_{\rm max}$. 
Unlike the Les Houches problem, where the nonsmooth optimal value is equal to the smoothed optimal value, that is $f^* = f_{\mu}^*=0$, for any $\mu $ as $\mu \to 0$, the same statement is not true for the Max Eigenvalue problem.  The smoothed Max Eigenvalue problem requires a complete eigendecomposition in order to obtain every eigenvalue for a given matrix, and since CVX does not allow such functions but only those that it knows to be convex such as the maximum eigenvalue function, we could not compute the optimal value $f_{\mu}^*$ from CVX. Instead, we use full BFGS with the max number of iterations set to $10^{5}$ to minimize \eqref{smax_eig} to high accuracy: we denote this computed value by $f_\mu^B$.

\begin{figure} 
	\centering
	\begin{tabular}{cc}
	\includegraphics[scale=0.3]{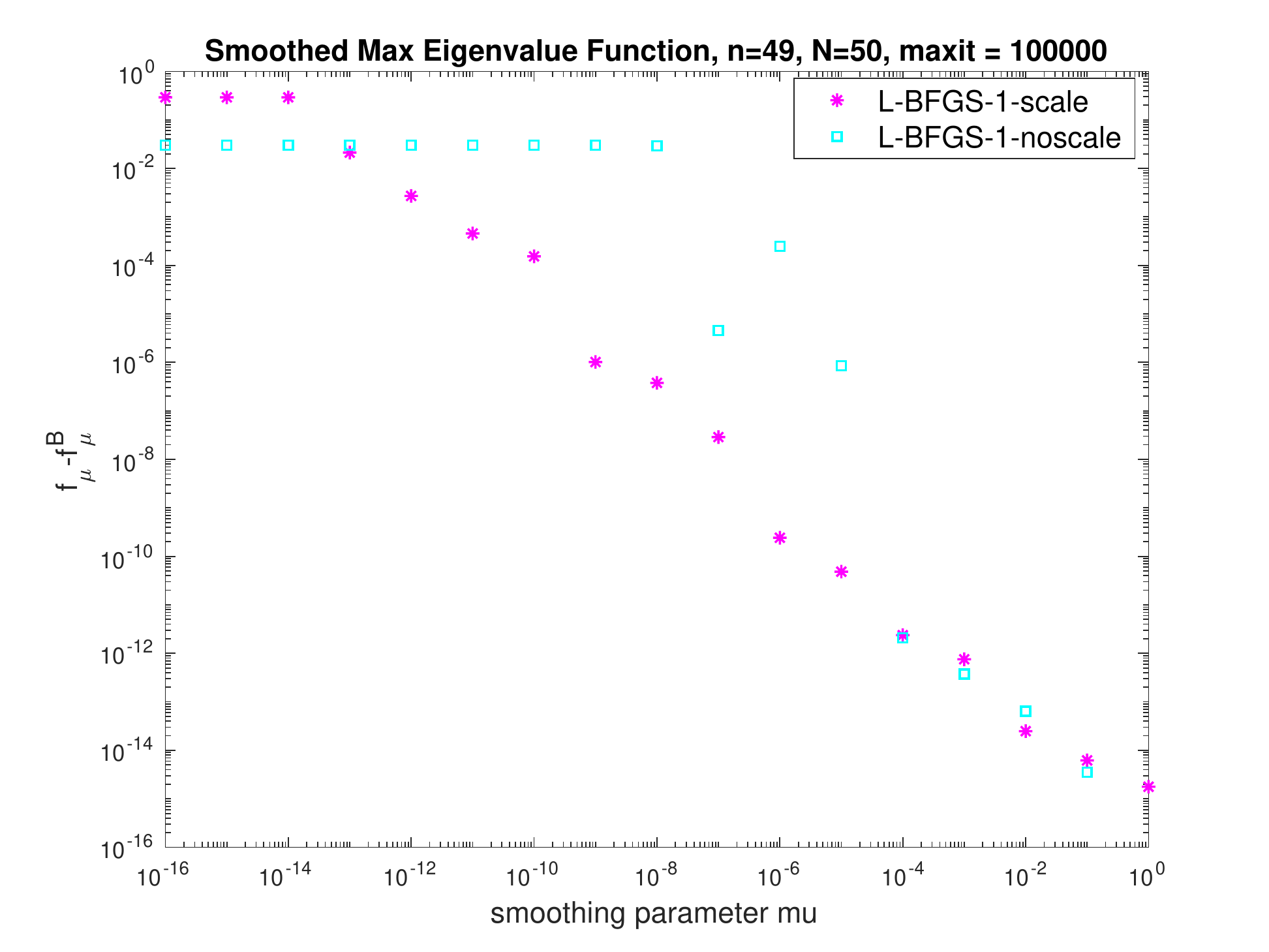}
	&
	 \includegraphics[scale=0.3]{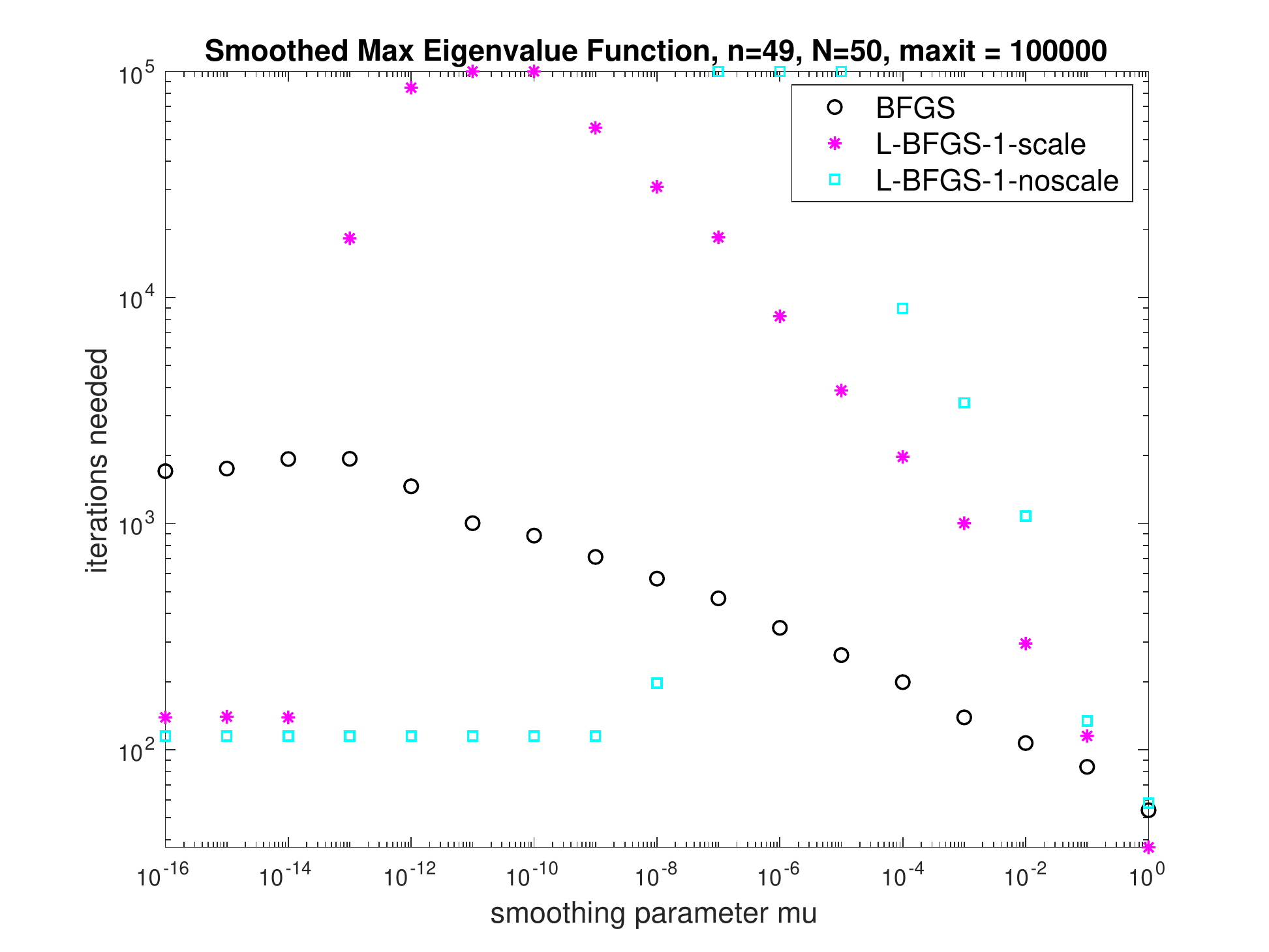} 
	\end{tabular}
	\caption
	[Smoothed Random Max Eigenvalue Problem -L-BFGS-1]
	{Comparing  L-BFGS-1 with and without scaling on the smoothed Max Eigenvalue problem \eqref{smax_eig} for  $N=50$ and  $n=49$. The left panel shows the final function value, shifted by $f^B_\mu$, the optimal value computed by BFGS, and the right panel shows the iteration count, both as a function of the smoothing parameter $\mu$. The maximum number of iterations is set to $10^5$.}
	\label{fig:sme}
\end{figure}

In Figure \ref{fig:sme}  we report on an experiment using the smoothed version of the same instance of the randomly generated Max Eigenvalue problem as earlier with $N=50$ and $n=49$, using L-BFGS-1 to minimize  \eqref{smax_eig}. The left panel shows the final value computed by scaled (magenta asterisks) and unscaled (cyan squares) L-BFGS-1 shifted by $f^B_{\mu}$ (the answer found by full BFGS), as a function of the smoothing parameter $\mu$, in log-log scale.
The right panel shows the number of iterations as a function of $\mu$, also in log-log scale.
The maximum number of iterations is $10^5$.

In the left panel we see that for $\mu=1$ down to $\mu=10^{-4}$ both methods yield about the same accuracy as each other, but that this deteriorates as $\mu$ decreases. Scaled L-BFGS-1 continues to obtain a reasonable approximation to the presumed accurate solution $f^B_\mu$ for $\mu$ down to $10^{-9}$, although this accuracy continues to decrease as $\mu$ is reduced.
Looking at the right panel, we see that starting with  $\mu=10^{-10}$ scaled L-BFGS-1 hits the maximum iteration limit and starting with $10^{-12}$ it breaks down before reaching the maximum iteration limit. In contrast,
unscaled L-BFGS-1 hits the maximum iteration limit for $\mu=10^{-5}$ and breaks down for $\mu\leq 10^{-8}$.

Table \ref{table:sme} shows the top 6 eigenvalues of the final answer found by BGFS and L-BFGS-1 for the smoothed Max Eigenvalue problem with $\mu=10^{-7}$.  We also repeat the optimal top 6 eigenvalues of the minimizer of the original nonsmooth function $f$ found by SDPT3 for the sake of comparison. Note that the result computed by applying scaled L-BFGS-1 to the smoothed problem agrees with the nonsmooth optimal value $f^*$ to 8 digits, compared to 0 digits when applied directly to the nonsmooth problem (Table \ref{table:mee}).
\begin{table} 
	\centering
	\footnotesize
	\begin{tabular}{|c|c|c|c|c|}
		\hline
		SDPT3 & BFGS  & Sc L-BFGS-1 & No L-BFGS-1 \\
		\hline
		7.82702970305352&7.82702976093363& 7.82702978971152& 7.82703432112405\\
		7.82702970305349&7.82702968960443& 7.82702971836333& 7.82703424950540\\
		7.82702970305348&7.82702966768912& 7.82702969644540& 7.82703422776180\\
		7.82702970305346&7.82702963829977& 7.82702966707913& 7.82703419808905\\
		7.82702970305334&7.82702954704101& 7.82702957585856& 7.82703410710019\\
		7.70350538019538&7.70350539089203& 7.70374015675041& 7.69886385782543\\
		\hline
	\end{tabular}
	\normalsize
	\caption
	[Top Eigenvalues of Smoothed Random Max Eigenvalue Problem -L-BFGS-1]
	{Top 6 eigenvalues of $C-\A^T y$ for Max Eigenvalue problem where $y$ is computed by applying  BFGS, scaled L-BFGS-1 and  unscaled L-BFGS-1 to $f_{\mu}$ with $\mu=10^{-7}$, for the same instance of the randomly generated Max Eigenvalue problem  as in Table \ref{table:mee}. The optimal multiplicity is 5. The first column gives the top 6 eigenvalues of the solution to the original nonsmooth problem.}
	\label{table:sme}
\end{table}

\begin{figure}
	\centering
	\begin{tabular}{cc}
	\includegraphics[scale=0.3]{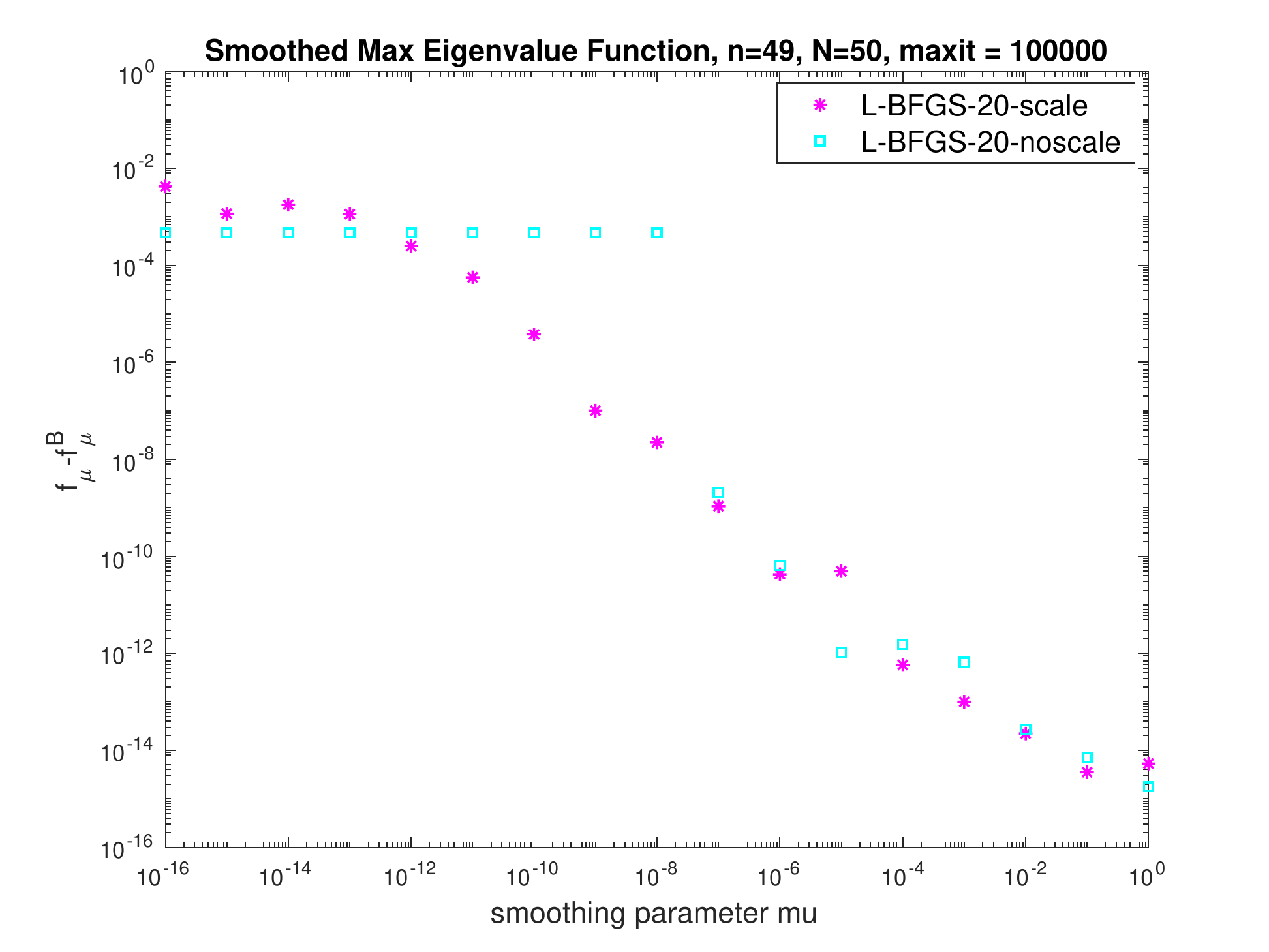}
	&
	\includegraphics[scale=0.3]{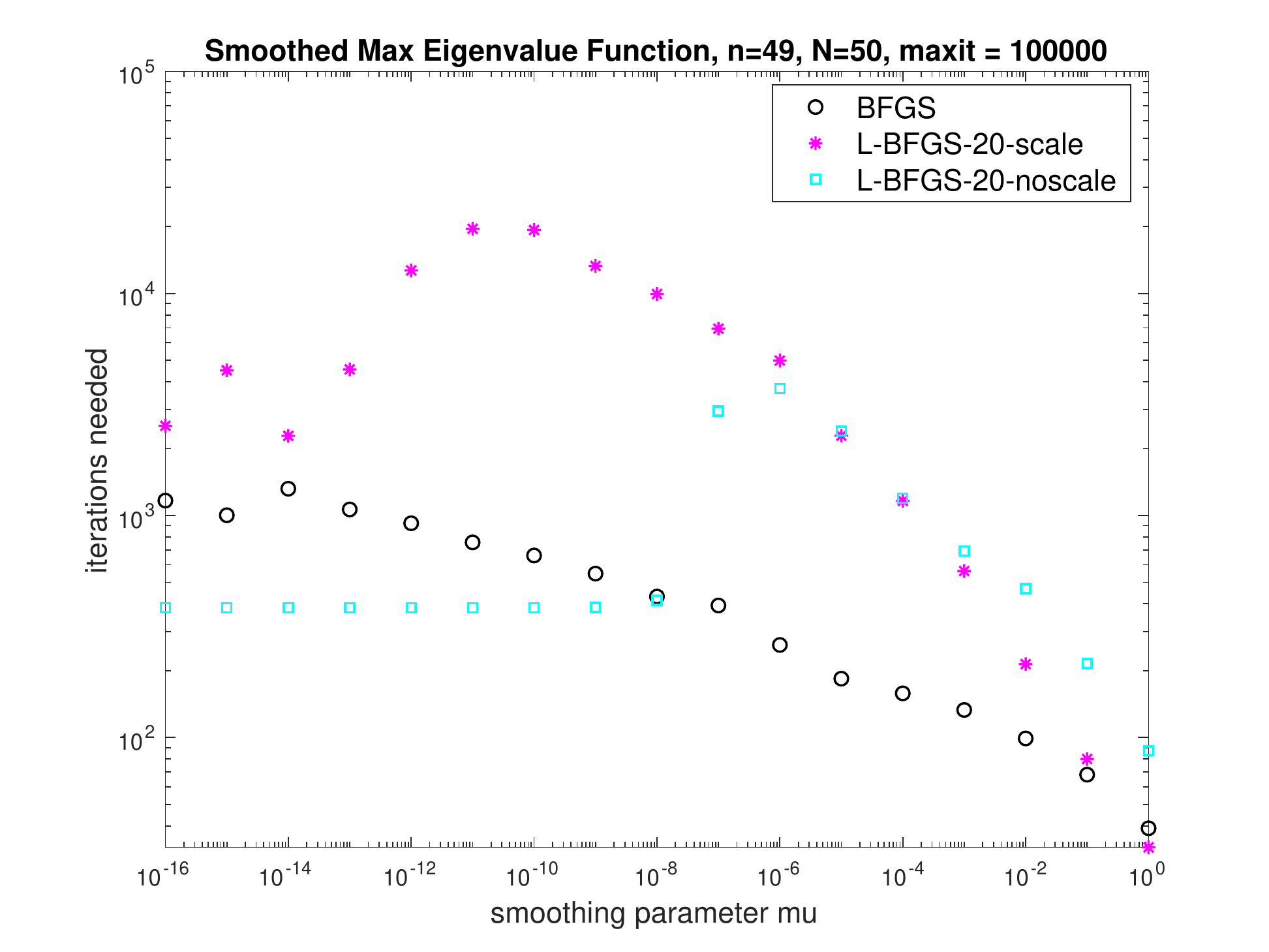}
	\end{tabular}
	\caption
	[Smoothed Random Max Eigenvalue Problem -L-BFGS-20]
	{Comparing  L-BFGS-20 with and without scaling on the smoothed Max Eigenvalue problem \eqref{smax_eig} for  $N=50$ and  $n=49$. The left panel shows the final function value, shifted by $f^B_\mu$, the optimal value computed by BFGS, and the right panel shows the iteration count, both as a function of the smoothing parameter $\mu$. The maximum number of iterations is set to $10^5$.}
	\label{fig:sme20}
\end{figure}

We repeated this experiment with $m=20$, reported in Figure \ref{fig:sme20}. In the left panel, we observe that 
the loss of accuracy in L-BFGS as a function of $\mu$ is less pronounced with 20 updates. Roughly speaking, overall the error decreases  by a factor of  $10^{-2}$.  In the right panel we see that neither scaled nor unscaled L-BFGS-20 reaches the maximum iteration limit, but that both methods break down for sufficiently small $\mu$. 
The top eigenvalues produced by L-BFGS-20 for $\mu=10^{-7}$ are shown in Table \ref{table:sme20}.

Comparing the final computed maximum eigenvalue in
Tables \ref{table:mee}, \ref{table:mee20}, \ref{table:sme}, and \ref{table:sme20}, note that full BFGS obtains a more accurate
solution when applied to the original nonsmooth problem than it does when applied to the smoothed approximation using $\mu=10^{-7}$, while the opposite is true for scaled L-BFGS-1 and scaled L-BFGS-20.
In summary, as with the Les Houches problem, it is much more effective to apply L-BFGS to the smoothed max eigenvalue problem than directly to the nonsmooth problem. As earlier, this is in sharp contrast to the behavior of full BFGS.

\begin{table}
	\centering
	\footnotesize
	\begin{tabular}{|c|c|c|c|c|}
		\hline
		Sc L-BFGS-20 & No L-BFGS-20 \\		
		\hline
		7.82702976201688 &7.82702976326908 \\
		7.82702969067290 &7.82702969273000 \\
		7.82702966877706 &7.82702966880138 \\
		7.82702963938214 &7.82702964120678 \\
		7.82702954813253 &7.82702954355970 \\
		7.70348252862186 &7.70344821498698 \\
		\hline
	\end{tabular}
	\normalsize
	\caption
	[Top Eigenvalues of Smoothed Random Max Eigenvalue Problem -L-BFGS-20]
	{Top 6 eigenvalues of $C-\A^T y$ for Max Eigenvalue problem where $y$ is computed by applying  BFGS, scaled L-BFGS-20 and  unscaled L-BFGS-20 to $f_{\mu}$ with $\mu=10^{-7}$, for the same instance of the randomly generated Max Eigenvalue problem  as in Table \ref{table:mee}. The optimal multiplicity is 5.}
	\label{table:sme20}
\end{table}

\subsection{Semidefinite Programming}

Consider the following primal and dual semidefinite programs (SDP) in standard form \cite{helRen, helOv}
\begin{align}
&\max_{X\in S^N} ~~~\langle	C,X \rangle  \label{sdpp}\\
&\mathrm{subject~to} ~~~~~ \A X=b  ~~~\mathrm{and} ~~~~ X\in S^N_+,\nonumber \\
&\min_{y \in \R^n} ~~~	b^Ty \label{sdpd} \\
&\mathrm{subject~to} ~~~ Z  = \A^Ty -C ~~~\mathrm{and} ~~~~{Z \in S^N_+},\nonumber
\end{align}
where $b \in \R^n$, $C \in S^N$ and
$\A:S^N\to \R^n$ is a linear operator as defined in \eqref{Alinopdef} and \eqref{Alinopdef2}.
Here $S^N_+\subseteq S^N$ denotes the cone of positive semidefinite $N\times N$ matrices.
Let us assume that strong duality holds, so that the optimal primal and dual values are the same, and that the optimal values are attained.  It follows that if  $X^*$ is an optimal solution to the primal problem \eqref{sdpp} and $Z^*$ is an optimal solution to the  dual problem \eqref{sdpd}, we have $X^*Z^* =0$. Further assume that $X^*$ is nonzero, and consequently, $Z^*$ has at least one eigenvalue equal to zero.
Then  the dual problem  \eqref{sdpd} is equivalent to the following unconstrained eigenvalue optimization problem 	
\beq\label{min_cdfef}
\min_{y \in \R^n} ~~~f(y),
\eeq
with  the exact penalty  dual function \cite{madel}
\beq \label{cdfef}
f(y) = b^Ty + \alpha  \max\{\lambda_{\rm max}(C - \A^Ty),~ 0\}, 
\eeq
for sufficiently large 
$\alpha$, where $\lambda_{\rm max}$ denotes maximum eigenvalue as earlier.
Note that this exact penalty function differs from the eigenvalue optimization formulation in \cite{helRen}, namely 
$$b^T y + \alpha \lambda_{\rm max}(C - \A^Ty)$$
which does not include the $\max\{\cdot~,~0\}$ operator. In that formulation, to give a correct equivalence $\alpha$ must be exactly equal to a critical value, as opposed to greater than or equal to this value. 
For the SDP problems we consider in the following subsections we already know valid lower bounds for $\alpha$. 
Note that at an optimal solution $y^*$  the maximum eigenvalue of $-Z^*=C - \A^T y^*$ is zero, often with multiplicity greater than one, and hence $f$ is nonsmooth at $y^*$.

\subsection{Max Cut Problem}\label{subsec:maxcut}
Our first example of semidefinite programming (SDP) arises from the Max Cut problem. This subsection and a subsequent one on the Matrix Completion problem were motivated by the recent paper \cite{madel} and the observation  made there that the first-order algorithms they used to minimize the penalized dual function \eqref{cdfef} arising from Max Cut  and Matrix Completion SDP relaxations were slow. Here we compare full BFGS, scaled and unscaled L-BFGS and the subgradient method (again with $t_k=1/k$) on penalized dual functions arising from 
Goemans-Williamson SDP relaxations of the Max Cut problem. We note that one of the key ideas in \cite{madel} is that, when the primal SDP optimal solution $X^*$ has rank much less than $N$, an accurate estimate of the optimal value of the SDP obtained from minimizing the penalized dual function allows the use of a novel method for obtaining efficient low-rank solutions to the \emph{primal} SDP even when $N$ is large.		

The primal  Max Cut SDP relaxation and its dual are \cite{GW95}: 
\begin{align}
&\max_X ~~~\frac{1}{4}\langle	L,X \rangle  \label{mxp}\\
&\mathrm{subject~to} ~~~~~ \diag X={\bf 1} ~~~\mathrm{and} ~~~~ X \in S^N_+ , \nonumber \\
&\min_{y \in \R^n} ~~~	{\bf 1}^Ty  \label{mxd}\\
&\mathrm{subject~to} ~~~ Z = \Diag y-\frac{1}{4}L  ~~~\mathrm{and} ~~~~ Z \in S^N_+ ,\nonumber
\end{align}
where $L$ is the Laplacian matrix of a given undirected graph, $\diag{}$ maps the diagonal
of a matrix to a vector, and $\Diag{}$ maps a vector to a diagonal matrix.
Note that these are instances of the primal and dual SDP introduced in \eqref{sdpp} and \eqref{sdpd}, respectively. By definition for the SDP Max Cut problem we have $n=N$.
The exact penalty dual function \eqref{cdfef} for the Max Cut SDP relaxation is%
\beq\label{pen_mxcut} 
f(y) ={\bf 1}^Ty + \alpha ~ \max\{\lambda_{\rm max}(L-\Diag y),~ 0\}.
\eeq
Due to the constant trace property of the primal Max Cut SDP \eqref{mxp}, the trace (nuclear) norm of the primal optimal  solution is known
to be $N$, and hence any solution $y^*$ to the penalized dual max cut  problem \eqref{pen_mxcut} with $\alpha \ge N$ is also a solution to the dual  SDP \eqref{mxd} and  vice versa \cite[Lem. 6.1]{madel}. 

We picked graph  $G1$ from the Gset group in the  sparse matrix collection \cite{Gset}  for the following experiment. $G1$ is  an unweighted  graph with $N=800$ vertices and the adjacency matrix is a sparse symmetric  matrix with $38352$ nonzero entries (all equal to 1).  
Since $N$ is relatively small we can apply the SDPT3 solver via CVX \cite{gb08} to the primal SDP \eqref{mxp}, obtaining the optimal primal and dual value $f^* = 12083.19765$. The rank $r^*$ of the optimal primal solution $X^*$ is $13$, and strict complementarity holds, so the nullity of the dual solution $Z^*$ is also $13$.

\begin{figure} 
	\includegraphics[width=5.5in,height=1.5in]{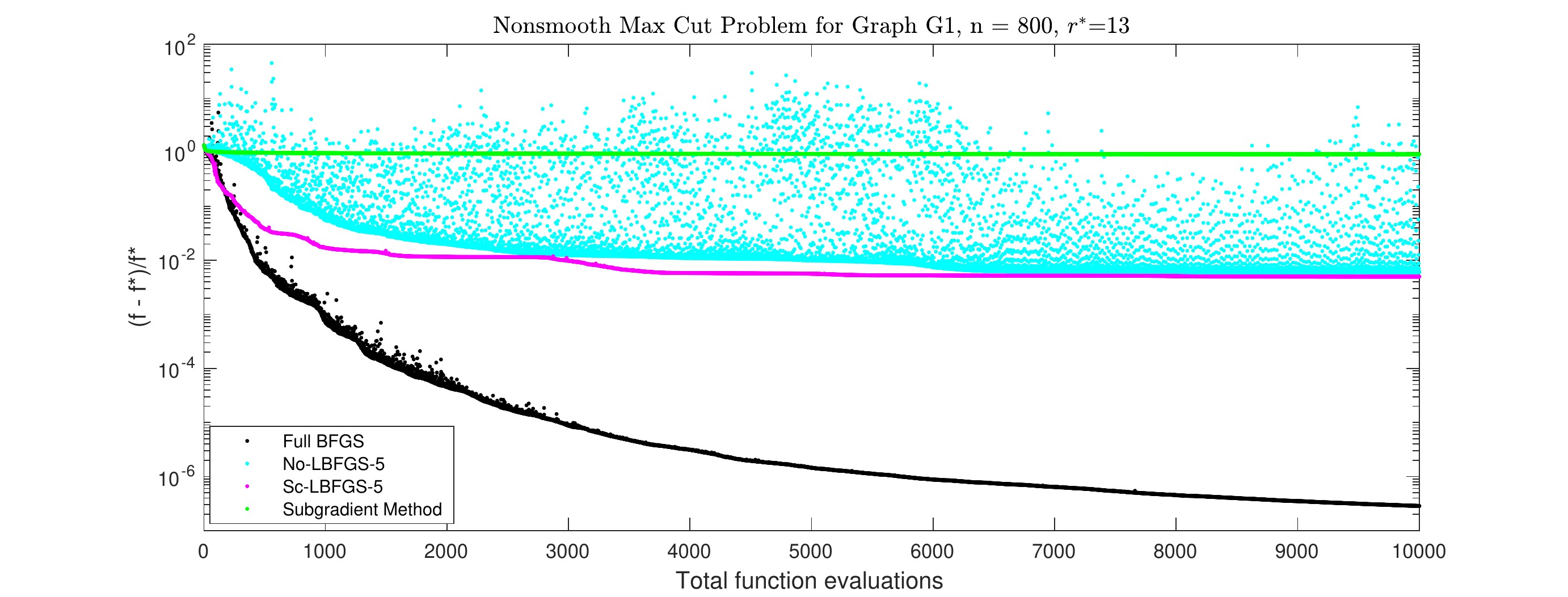}
	\caption
	[Penalized Dual Max Cut Problem -L-BFGS-5]
	{Comparing BFGS, L-BFGS-5 with and without scaling and the subgradient method on the penalized dual Max Cut  problem \eqref{pen_mxcut}.}
	\label{fig:mx1}
\end{figure}

In  Figure \ref{fig:mx1}, we compare the performance of full BFGS, L-BFGS-5 with and without scaling, and the subgradient method (with $t_k=1/k$) to minimize the penalized dual function \eqref{pen_mxcut} with $\alpha  = 2N=1600$. We display all the function values that were computed, with
the maximum number of function evaluations set to $10^4$.
The vertical axis shows the relative error $(f-f^*)/f^*$. 
In  contrast to the two experiments 
presented in Figure \ref{fig:me1} and \ref{fig:me20} for the nonsmooth Max Eigenvalue problem, the results of this experiment are in favor of L-BFGS when compared to the subgradient method.
Both scaled L-BFGS-5 (magenta dots) and unscaled L-BFGS-5 (cyan dots) reduce the relative error down to  
below $10^{-2}$ while the subgradient method (blue dots) reduces the relative error to about $10^{0}$. 
Full BFGS (black dots) reduces the relative error to below $10^{-6}$. 

\begin{figure} 
        \includegraphics[width=5.5in,height=1.5in]{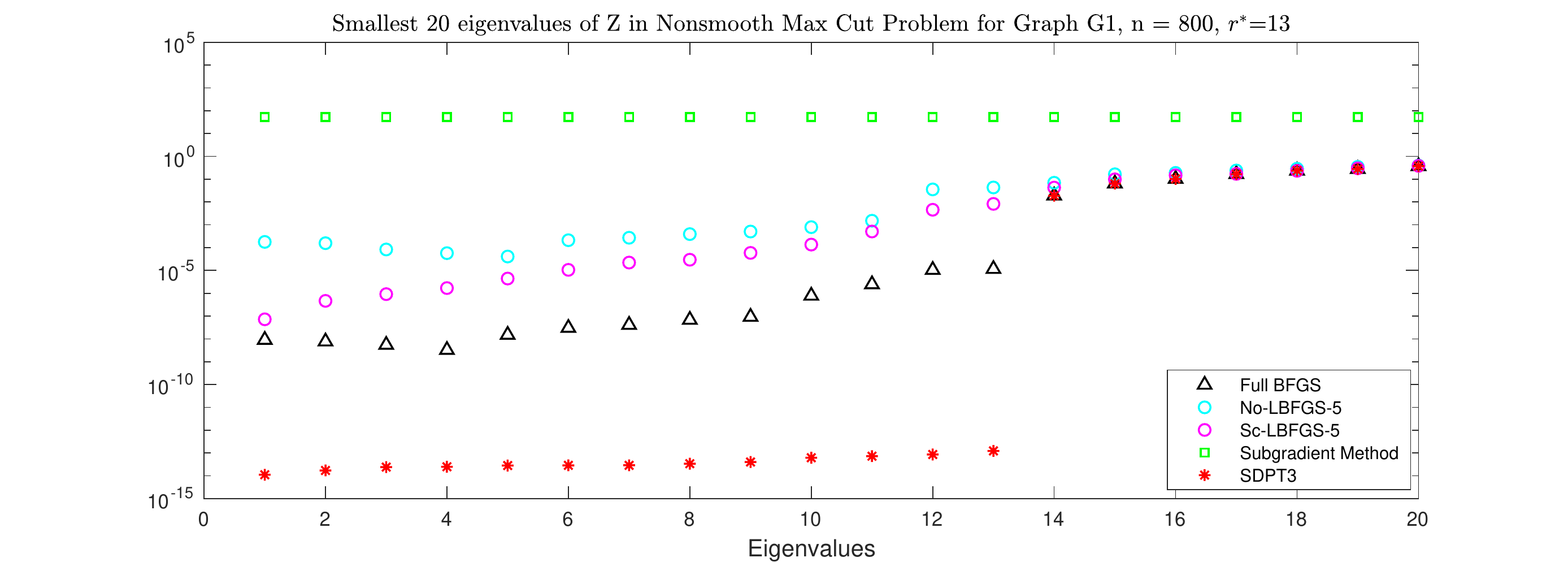}
	\caption
	[Penalized Dual Max Cut Problem -Eigenvalues of the Dual Slack Matrix -L-BFGS-5]
	{Comparing smallest 20 eigenvalues of the dual slack matrix $Z$ obtained by BFGS, L-BFGS-5 with and without scaling and the subgradient method on the penalized dual Max Cut  problem \eqref{pen_mxcut} for  the $G1$ graph with $n=800$.  The nullity of the optimal dual slack matrix $Z^*$ is 13. The smallest 20 eigenvalues obtained from SDPT3 are shown as well. The lack of monotonicity at the left end of some of the plots occurs because we actually plotted the absolute values of the ordered largest eigenvalues of $-Z$, and some of these eigenvalues are positive, either because of rounding errors or insufficient accuracy in the optimization. } 
	\label{fig:nmxi}
\end{figure}

In Figure \ref{fig:nmxi}, we show the \emph{negative} of the top 20 eigenvalues of the final negative dual slack matrix $-Z$ (equivalently, the smallest 20 eigenvalues of $Z$) obtained by the four methods, along with values obtained by SDPT3. It is interesting to note that full BFGS approximates the eigenvalues well, in the sense that it clearly separates the first 13 approximately zero eigenvalues from the approximations to the nonzero eigenvalues, although not as decisively as SDPT3. However, the final eigenvalues obtained by L-BFGS-5 are not clearly separated. 

\begin{figure} 
	\includegraphics[width=5.5in,height=1.5in]{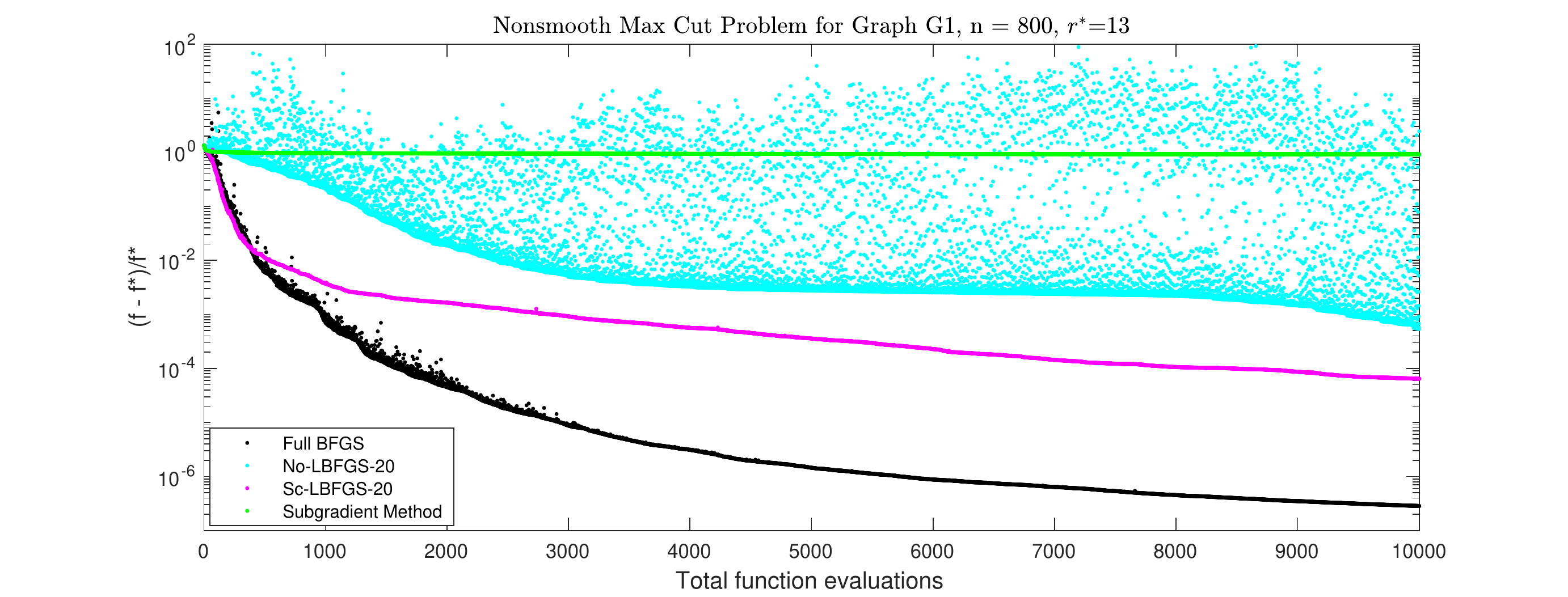}
	\caption
	[Penalized Dual Max Cut Problem -L-BFGS-20]
	{Comparing BFGS, L-BFGS-20 with and without scaling and the subgradient method on the penalized dual Max Cut  problem \eqref{pen_mxcut}. }
	\label{fig:mx2}
\end{figure}

\begin{figure} 
	\includegraphics[width=5.5in,height=1.5in]{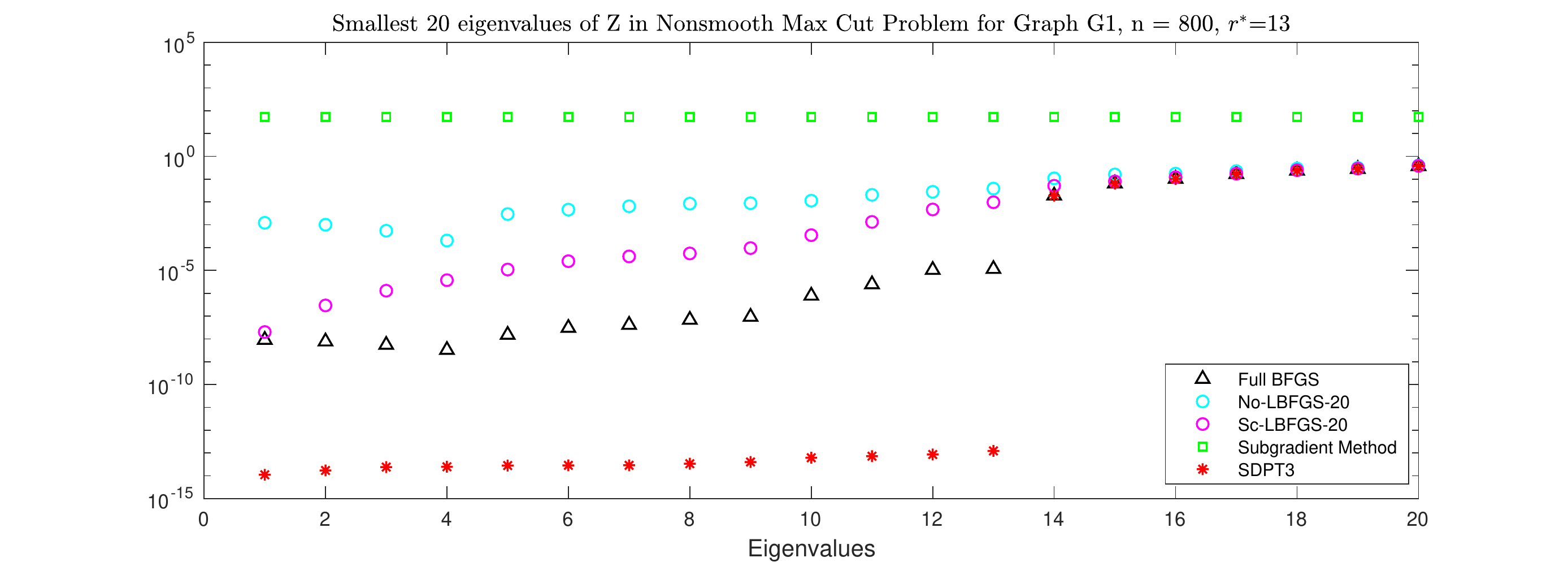}
	\caption
	[Penalized Dual Max Cut Problem -Eigenvalues of the Dual Slack Matrix -L-BFGS-20]
	{Comparing smallest 20 eigenvalues of the dual slack matrix $Z$ obtained by BFGS, L-BFGS-5 with and without scaling and the subgradient method on the penalized dual Max Cut  problem \eqref{pen_mxcut} for  $G1$ graph with $n=800$. The nullity of the optimal dual slack matrix $Z^*$ is 13. The smallest 20 eigenvalues obtained from SDPT3 are shown as well. See legend of Figure \ref{fig:nmxi} regarding eigenvalue monotonicity.}
	\label{fig:nmxi2}
\end{figure}

Next we increase the number of L-BFGS updates from $m=5$ to $m=20$, showing the results in Figures \ref{fig:mx2} and \ref{fig:nmxi2}. The results for BFGS and the subgradient method are shown again for comparison.
We see that scaled L-BFGS-20 now reduces the relative error down to  $10^{-4}$ and unscaled to about $10^{-3}$, compared to about $10^{-2}$ for scaled and unscaled L-BFGS-5.
However, the eigenvalue plot is similar to the corresponding plot for L-BFGS-5:  neither variant is able to discover that the nullity of the optimal dual slack matrix $Z^*$ is 13.

\subsection{Smoothed Max Cut Problem}\label{subsec:smax_cut}
Consider now Nesterov smoothing of the penalized dual Max Cut  problem \eqref{pen_mxcut}:
\beq\label{smax_cut}
f_{\mu}(y) = {\bf 1}^Ty +  \alpha \mu \log \left( 1 + \sum_{i=1}^n  \exp \left( \lambda_i(L-\Diag y)/\mu \right)  \right )
- \alpha \mu\log (n+1).
\eeq
Note the presence of the term ``1" which does not appear in \eqref{smax_eig}: this reflects the presence of the $\max\{\cdot~,~0\}$ operator in the penalty function \eqref{cdfef}.

\begin{figure} 
	\includegraphics[width=5.5in,height=1.5in]{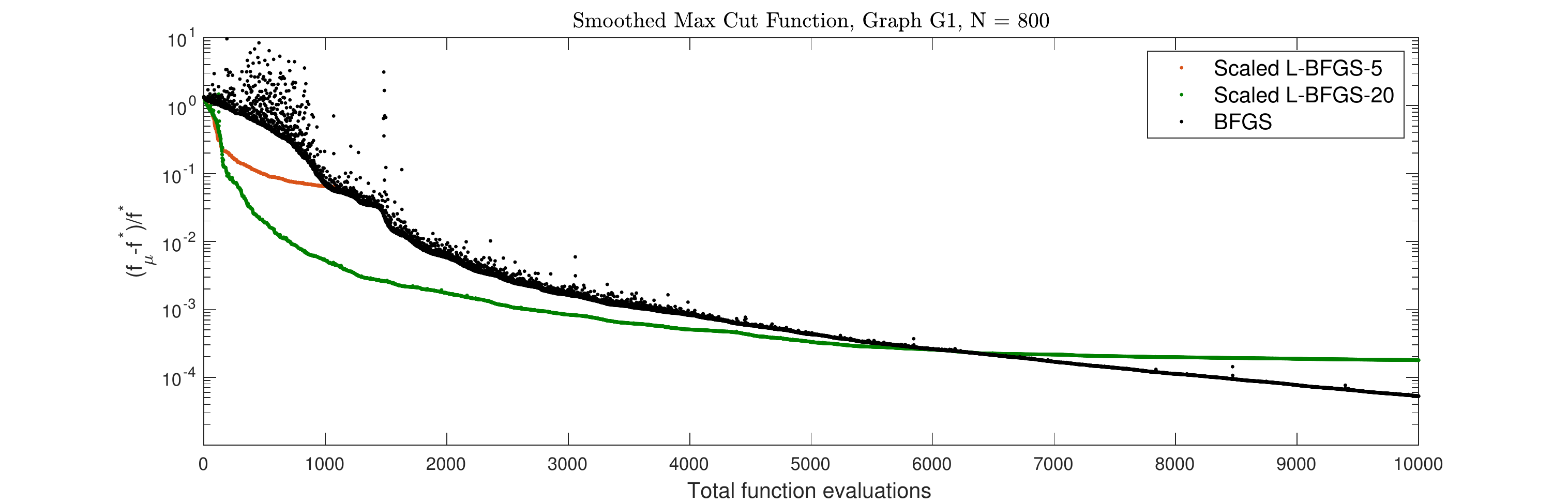}
	\caption
	[Smoothed Max Cut Problem]
	{Comparing scaled L-BFGS-5,  scaled L-BFGS-20 and BFGS on Smoothed Max Cut problem \eqref{smax_cut} for the same problem as in the Nonsmooth Max Cut problem. The smoothing parameter is $\mu=10^{-7}$. }
	\label{fig:s_maxcut}
\end{figure}

Figure \ref{fig:s_maxcut} shows the results of applying BFGS, scaled L-BFGS-5 and scaled L-BFGS-20 to minimize \eqref{smax_cut} with $\mu=10^{-7}$.
We do not include the unscaled L-BFGS variants because it seems clear that they offer no advantage on smooth problems. 
Interestingly, and in contrast to our other experiments with smoothed nonsmooth functions, scaled L-BFGS-20 does better than
full BFGS during its first several thousand function evaluations, but eventually it is overtaken by full BFGS. Scaled L-BFGS-5 does relatively poorly. 
\begin{table} 
\centering
\scriptsize
\begin{tabular}{ |c|c|}
\hline
SDP-optimal   &   12083.19765454945     \\
\hline
BFGS-smooth &       12083.83341694415 \\
\hline
BFGS-nonsmooth &     12083.20108505506  \\
\hline
L-BFGS-20-smooth &    12085.35533081401   \\
\hline
L-BFGS-20-nonsmooth &    12083.97779371002  \\
\hline
L-BFGS-5-smooth &      12848.47893036591 \\
\hline
L-BFGS-5-nonsmooth &    12143.81352524515  \\
\hline
\end{tabular}
\caption
[Max Cut Problem - Objective Values]
{Final objective value we obtained from   full BFGS, scaled L-BFGS-20 and scaled L-BFGS-5 on the  smoothed and nonsmooth Max Cut problem presented in Figures \ref{fig:s_maxcut},  \ref{fig:mx2} and \ref{fig:mx1} respectively. The optimal SDP value $f^*$ is also shown for comparison.}
\label{table:objval0}
\end{table}

Table \ref{table:objval0} shows the final answers we obtained from full BFGS, scaled L-BFGS-20 and scaled L-BFGS-5 on the smoothed and nonsmooth Max Cut problem. The optimal SDP value $f^*$ is also shown.
All three of full BFGS, scaled L-BFGS-5 and scaled L-BFGS-20 obtain better approximations to $f^*$
when applied directly to the nonsmooth problem than when applied to the smoothed problem with $\mu=10^{-7}$.
This is in contrast to what we observed for the Les Houches problem and the Max Eigenvalue problem,
where this was true only for full BFGS.

It would be interesting to investigate whether L-BFGS might be useful in the solution of large-scale Max Cut problems or other
SDPs. A first step in this direction appears in \cite[Sec.\ 4.2.5]{AA20}, utilizing a variant of the smoothed
exact penalty dual function \eqref{smax_cut} that includes only the largest eigenvalues in the smoothing, 
since these are the ones that dominate \eqref{smax_cut}. This allows the use of
\matlab's {\tt eigs} to compute only the largest few eigenvalues of $C-\A^T y$ via the Lanczos method.


\subsection{Matrix Completion Problem}\label{subsec:matcom}

The Matrix Completion problem is as follows. Suppose $\X \in \R^{N_1\times N_2}$ denotes a low-rank matrix for which we only have access to  some of its entries and would like to recover entirely by minimizing the rank over all matrices whose entries agree with the known values. This rank minimization problem is NP-hard, so we relax it by minimizing a well-known convex surrogate for the rank: the nuclear norm (sum of all the singular values). Let $\Omega$ be  the set of pairs $(i,j)$ for which $\X_{ij}$ is known. 
Then the nuclear norm minimization problem can be expressed as the following SDP \cite{FazRech}
\begin{align}
&\max_{X \in S^{(N_1 + N_2)} }~~~-\tr{W_1}- \tr {W_2}  \label{mcp}\\
&\mathrm{subject~to} ~~~  U_{ij}= \X_{ij},~(i,j)\in \Omega,\nonumber \\
& X = \left[
\begin{array}{cc}
W_1 & U\\
U^T & W_2
\end{array}
\right] \in S^{(N_1 + N_2)}_+. \nonumber
\end{align}
We write the primal problem in the \emph{max} form in order to   
be consistent with the SDP form \eqref{sdpp}, with $N = N_1+N_2$.  
Define the constraint $U_{ij}= \X_{ij}$ for $(i,j)\in \Omega$ in linear operator form $\mathcal B(U)=b$, where  
$b \in \R^n$ with $n=|\Omega|$ is the vector consisting of the known entries of $\X$ in some prescribed order and $\mathcal B: \R^{N_1\times N_2} \to \R^n$, with $\mathcal B^T$ its adjoint operator. The dual SDP is
\begin{align}
&\min_{y \in \R^n} ~~~	b^Ty  \label{mcd}\\
&\mathrm{subject~to} ~~~ Z= 	\left[
\begin{array}{cc}
I_{N_1} & \mathcal B^T(y)\\
\left ( \mathcal B^T(y)\right) ^T & I_{N_2}
\end{array}
\right],~~  Z \in S^{(N_1 + N_2)}_+ \green{.}	\nonumber
\end{align}
The exact penalty dual function \eqref{cdfef} 
for the Matrix Completion  problem is then 
\beq \label{pendmt}
f(y) =b^Ty + \alpha  \max \Bigg \{ \lambda_{\rm max} \left (- 	
\left[
\begin{array}{cc}
	I_{N_1} & \mathcal B^T(y)\\
	\left ( \mathcal B^T(y)\right) ^T & I_{N_2}
\end{array}
\right]
\right )  ,~ 0 \bigg \}.
\eeq

For the experiment in this part, we generated a low-rank 
random  matrix $\X$ of size $N_1=20$ by $N_2=160$ with rank $3$. 
We then selected the ordered pairs in $\Omega$ randomly with the probability of each $(i,j)$ being included set to $0.2$.
For this problem instance we got $|\Omega|=n=587$. 
We then applied the various methods to minimize  \eqref{pendmt} with $\alpha = 2\|X^*\|_*= 2\tr{X^*} = -2f^* =  3.0796$, where $f^* = -1.5398$ and $X^* \in S^{N_1 + N_2}_+ $ were obtained from solving the  SDP \eqref{mcd} via SDPT3. Note that the optimal value is negative because of the minus sign in the \emph{max} formulation of the primal SDP.

\begin{figure} 
	\includegraphics[width=5.5in,height=1.5in]{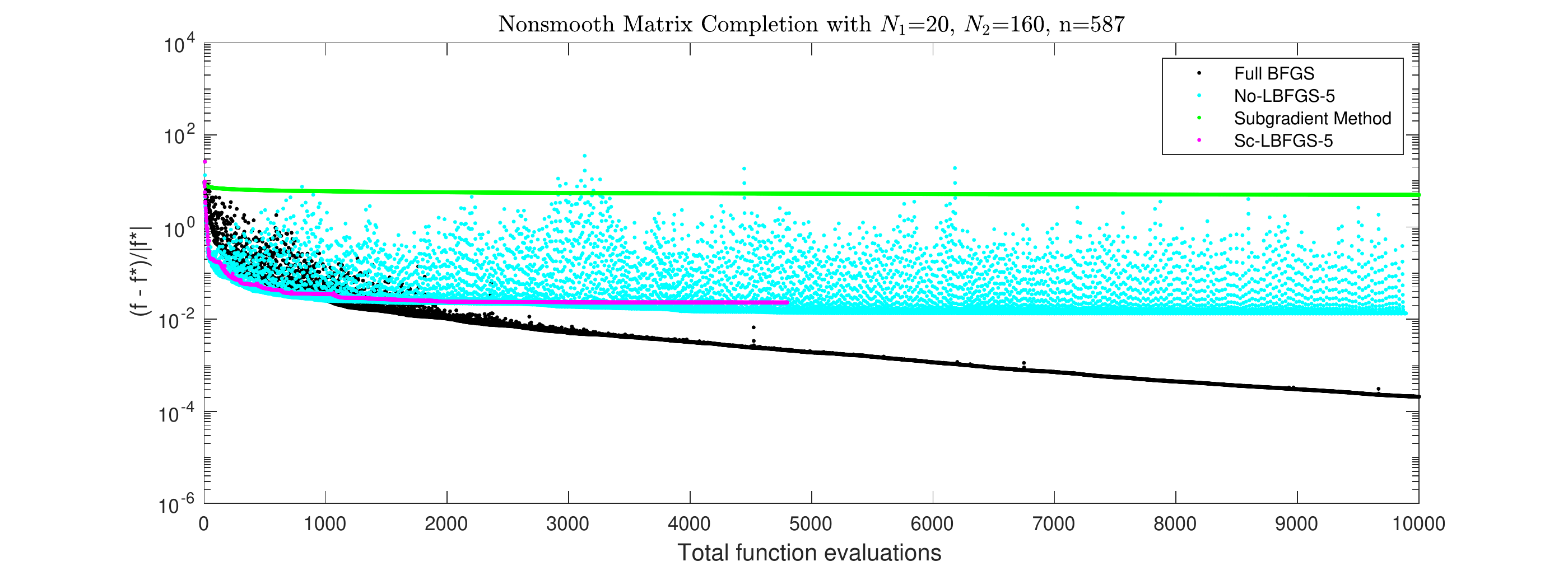}
	\caption
	[Penalized Dual Matrix Completion Problem -L-BFGS-5] 
	{Comparing BFGS, LBFGS-5 with and without scaling and the subgradient method on the penalized dual Matrix Completion problem \eqref{pendmt}.}
	\label{fig:mt1}
\end{figure}

Figure \ref{fig:mt1}  shows  the performance of full BFGS, L-BFGS-5 with and without scaling, and the subgradient method (with $t_k=1/k$). The vertical axis shows the relative error $(f-f^*)/|f^*|$. 
The maximum number of function evaluations is set to $10^4$. As is evident from the plot, both  variants of L-BFGS-5 outperform the subgradient method, even though scaled L-BFGS-5 quits early before 5000 evaluations  and unscaled L-BFGS-5 just before 10000 evaluations.

\begin{figure} 
	\includegraphics[width=5.5in,height=1.5in]{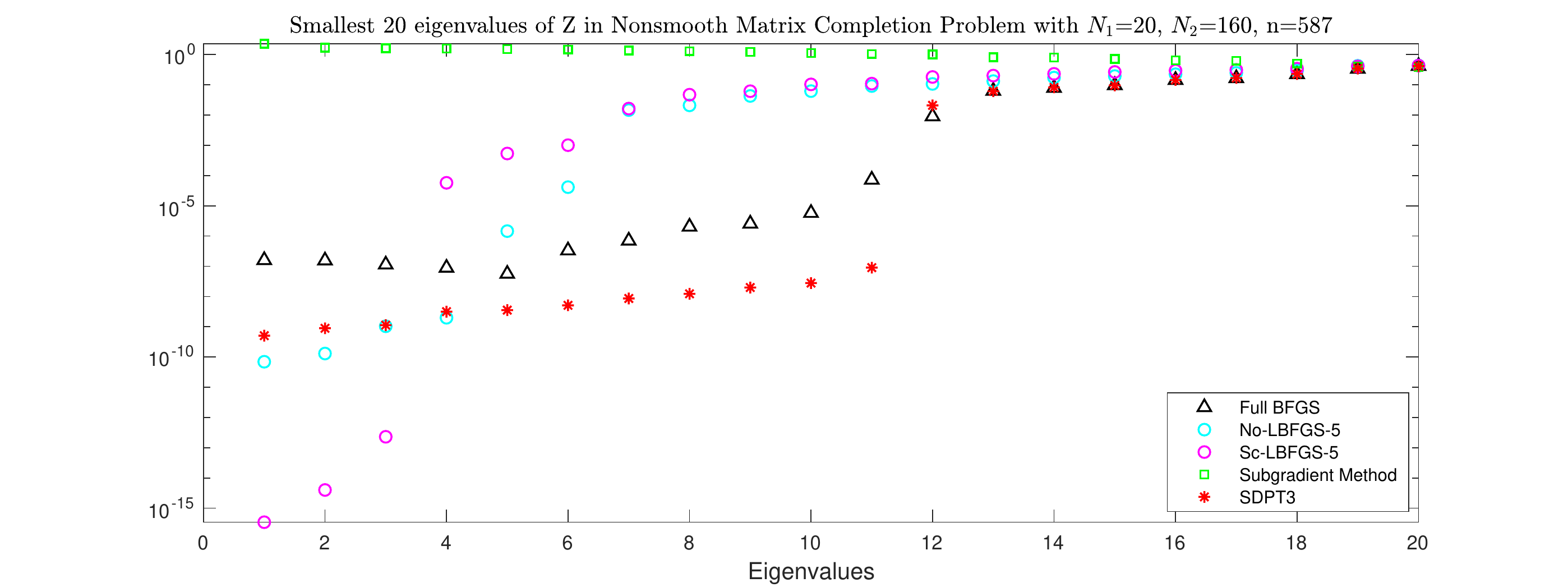}
	\caption
	[Penalized Dual Matrix Completion Problem -Eigenvalues  of the Dual Slack Matrix -L-BFGS-5] 
	{Comparing smallest 20 eigenvalues of the dual slack matrix $Z$ obtained by BFGS, L-BFGS-5 with and without scaling and the subgradient method on the penalized dual Matrix Completion problem \eqref{pendmt}. See legend of Figure \ref{fig:nmxi} regarding eigenvalue monotonicity.}
	\label{fig:mt1i}
\end{figure}

Figure  \ref{fig:mt1i} presents the negative of the top 20 eigenvalues of the final negative dual slack matrix $-Z$, or equivalently, the  20 smallest eigenvalues of $Z$, obtained by the four methods, along with values obtained by SDPT3. 
As before, BFGS is able to separate the zero and nonzero eigenvalues of $Z^*$, agreeing with SDPT3 that the nullity of $Z^*$ is effectively 11. Note that this is larger than $3$, the rank of the original matrix $\X$, implying that 20\% was not enough observations to reconstruct $\X$. It is interesting that the eigenvalues of the solution found by scaled L-BFGS-5 do suggest a nullity of 3, but this may just be a coincidence.

\begin{figure} 
	\includegraphics[width=5.5in,height=1.5in]{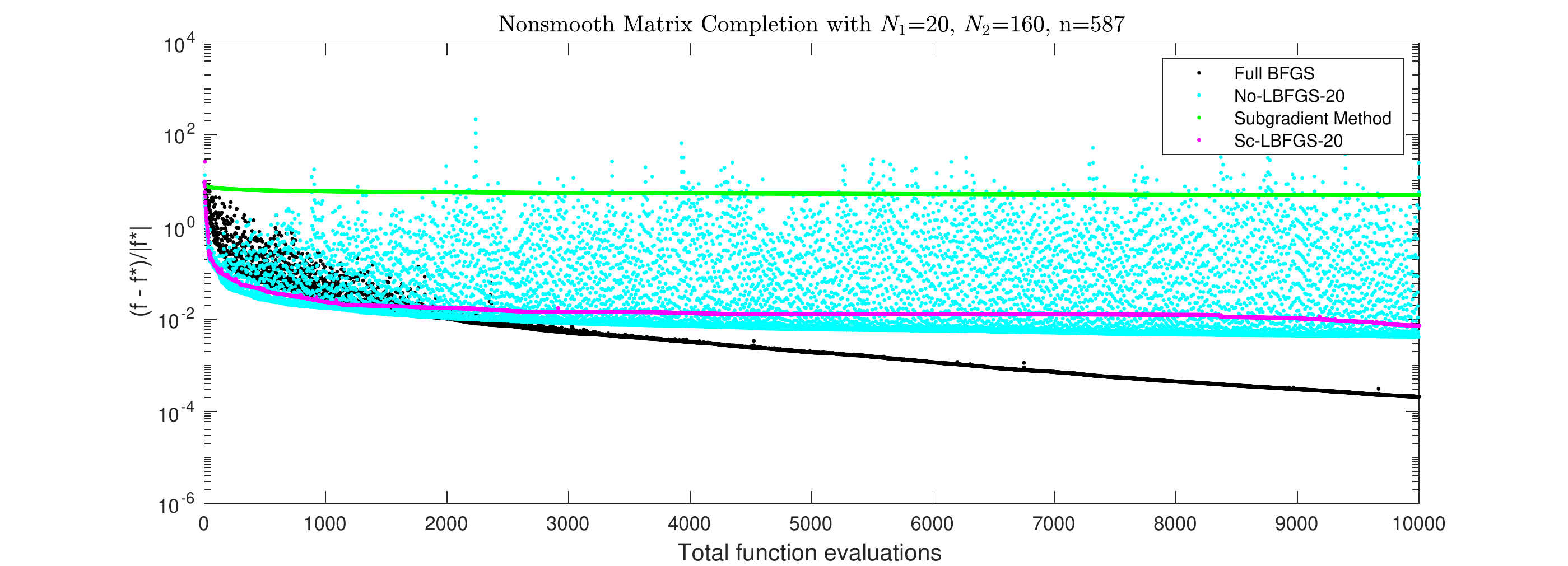}
	\caption
	[Penalized Dual Matrix Completion Problem -L-BFGS-20] 
	{Comparing BFGS, LBFGS-20 with and without scaling and subgradient method on the penalized dual Matrix Completion problem \eqref{pendmt}.}
	\label{fig:mt20}
\end{figure}

\begin{figure} 
	\includegraphics[width=5.5in,height=1.5in]{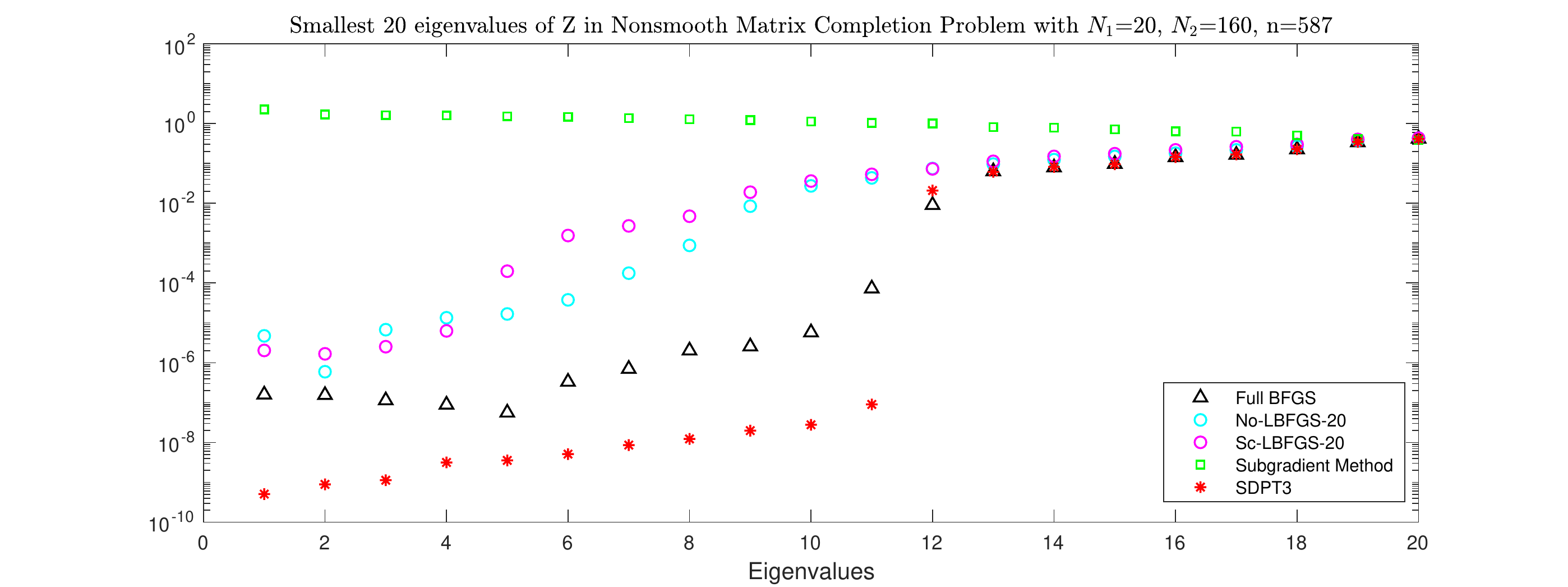}
	\caption
	[Penalized Dual Matrix Completion Problem -Eigenvalues  of the Dual Slack Matrix -L-BFGS-20] 
	{Comparing smallest 20 eigenvalues of the dual slack matrix $Z$ obtained by BFGS, L-BFGS-20 with and without scaling and the subgradient method on the penalized dual Matrix Completion problem \eqref{pendmt}. See legend of Figure \ref{fig:nmxi} regarding eigenvalue monotonicity.}
	\label{fig:mt20i}
\end{figure}

In the experiment reported in Figures  \ref{fig:mt20}  and \ref{fig:mt20i}, we increase $m$ to 20 and again we compare the relative error and the smallest 20 eigenvalues of the dual slack 
matrix, respectively. In both plots the result from BFGS and the subgradient method are  repeated for comparison. In Figure \ref{fig:mt20}, neither L-BFGS-20 method quits early this time; the unscaled variant gets a slightly lower answer. 
The eigenvalues shown for L-BFGS-20 in Figure \ref{fig:mt20i} do not suggest any conclusion about the nullity of $Z^*$.




\subsection{Smoothed Matrix Completion Problem}\label{subsec:smax_com}
Nesterov smoothing of the penalized dual Matrix Completion  problem \eqref{pendmt} gives 
\begin{align} 
&f_{\mu}(y) = b^Ty  \label{smat_com}\\
& + \alpha \mu\log \left( 1+ \sum_{i=1}^{N_1+N_2} \exp \left(\lambda_i \left(
- 	
\left[
\begin{array}{cc}
I_{N_1} & \mathcal B^T(y)\\   \nonumber
\left ( \mathcal B^T(y)\right) ^T & I_{N_2}
\end{array}
\right]
\right)
/\mu \right)  \right) \\
&- \alpha\mu\log(N_1+N_2 + 1).   \nonumber
\end{align}

\begin{figure} 
	\includegraphics[width=5.5in,height=1.5in]{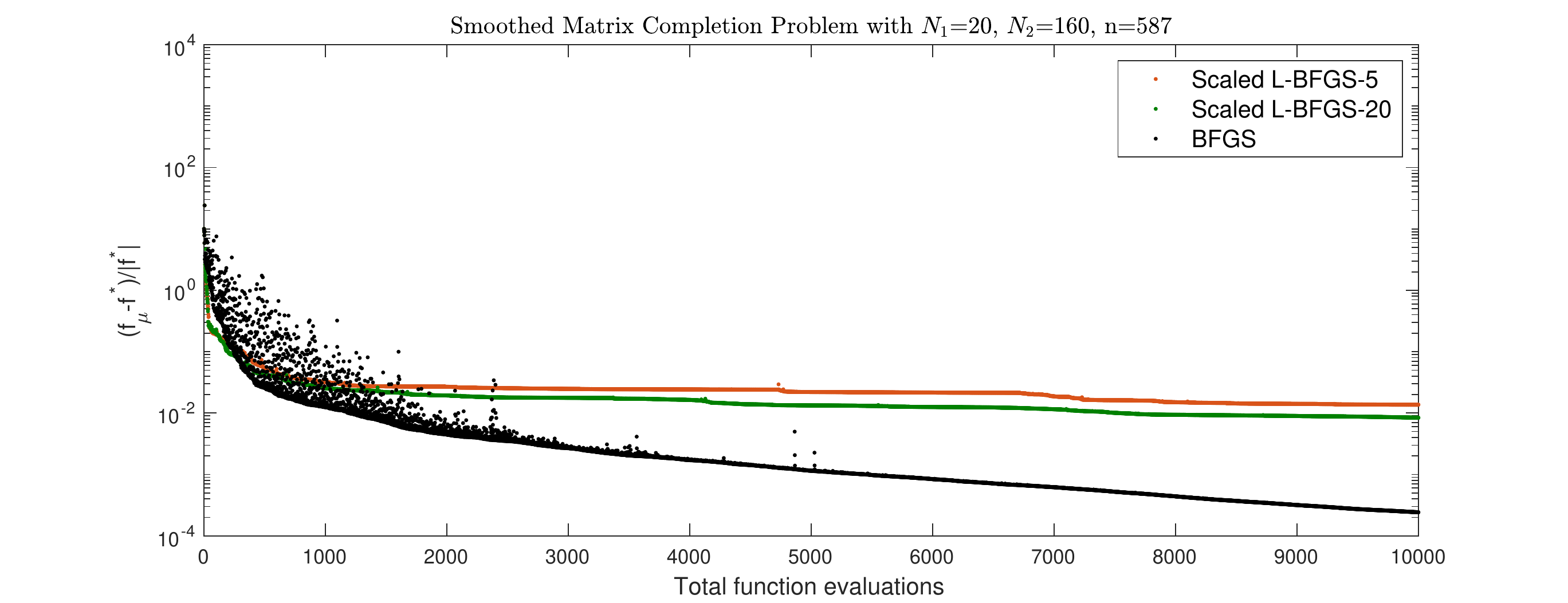}
	\caption
	[Smoothed Matrix Completion Problem]
	{Comparing scaled L-BFGS-5,  scaled L-BFGS-20 and BFGS on Smoothed Matrix Completion problem \eqref{smat_com} for the same problem as in the Nonsmooth Matrix Completion problem. The smoothing parameter is $\mu=10^{-7}$. }
	\label{fig:smc}
\end{figure}

In Figure \ref{fig:smc}, we show the result of applying full BFGS, scaled L-BFGS-5 and scaled L-BFGS-20 to \eqref{smat_com} with $\mu=10^{-7}$. The underlying reference matrix, $\X$, is the same matrix as in the nonsmooth experiment in Figure \ref{fig:mt1}, with $N_1=20$, $N_2=160$, $f^* = -1.5398$. We set $ \alpha = 3.0796$ as before.  The maximum number of function evaluations is set to $10^4$.  
 
We see in Figure \ref{fig:smc} that although
the relative error obtained by scaled L-BFGS-20 and full BFGS are about the same as when applied to the nonsmooth function, scaled L-BFGS-5 gets a lower error when it is applied to the smoothed function, and more importantly, it does not break down early. 

Table \ref{table:objval}
shows the final answers we obtained from   full BFGS, scaled L-BFGS-20 and scaled L-BFGS-5 on the  smoothed and nonsmooth Matrix Completion problem. The optimal SDP value $f^*$ is also shown for comparison.
\begin{table} 
	\centering
	\footnotesize
	\begin{tabular}{|c|c|}		\hline
		SDP-optimal   &   $-$1.53978555575175     \\
		\hline
		BFGS-smooth &       $-$1.53941270679104\\		
		\hline
		BFGS-nonsmooth		 &   $-$1.53946718487809    \\		
		\hline
		L-BFGS-20-smooth		 &    $-$1.52686081626082   \\		
		\hline
		L-BFGS-20-nonsmooth	 &    $-$1.52864965852708   \\		
		\hline				 
		L-BFGS-5-smooth &      $-$1.51868588317043 \\		
		\hline				 
		L-BFGS-5-nonsmooth &   $-1$.50422184883968   \\		
		\hline
	\end{tabular}
	\normalsize
	\caption
	[Matrix Completion Problem - Objective Values]
	{Final objective value we obtained from   full BFGS, scaled L-BFGS-20 and scaled L-BFGS-5 on the  smoothed and nonsmooth Matrix Completion problem presented in Figures \ref{fig:smc},  \ref{fig:mt20} and \ref{fig:mt1} respectively. The optimal SDP value $f^*$ is also shown for comparison.}
	\label{table:objval}
\end{table}

\section{Concluding Remarks}\label{sec:conclu}

In \S \ref{sec:theory}, we presented theoretical results showing that L-BFGS may converge to non-optimal points when applied to a simple
class of nonsmooth functions. In \S \ref{sec:practice}, we investigated whether the same phenomenon holds in practice. 
We found that when applied to a nonsmooth function directly,  L-BFGS, especially its scaled variant, often breaks down with a poor
approximation to an optimal solution, in sharp contrast to full BFGS. Unscaled L-BFGS is less prone to breakdown but conducts far more function evaluations per iteration than scaled L-BFGS does, and thus it is slow. Nonetheless, it is often the case that both variants obtain better results than the provably convergent, but slow, subgradient method.

On the other hand, when applied to a smooth approximation of a nonsmooth function, 
scaled L-BFGS invariably performs better than unscaled L-BFGS,
often obtaining good results even when the problem is quite ill-conditioned. In particular, scaled L-BFGS 
may be a reasonable approach to finding approximate minimizers of
smoothed exact penalty dual functions arising in large-scale semidefinite programs, although further investigation is needed to investigate the practicality of this approach. Minimization of the SDP exact penalty dual function is a key component of a recently proposed method for solving large-scale SDPs with low-rank primal solutions \cite{madel}.

Most importantly, we find that although L-BFGS is often a reliable method for minimizing ill-conditioned smooth problems, it frequently fails when the condition number is so large that the function is effectively nonsmooth. This behavior is in sharp contrast to the behavior of full BFGS, which is consistently reliable for nonsmooth optimization problems. We arrive at the conclusion that, for large-scale nonsmooth optimization problems for which full BFGS and other methods are not practical, it is 
often  better to apply L-BFGS to a smoothed variant of a nonsmooth problem than to apply it directly to the nonsmooth problem.

{\bf Acknowledgments.} The work of the first author was conducted as part of her Ph.D.\ studies at the Courant Institute of Mathematical
Studies, New York University. Many thanks to Margaret H.\ Wright for providing her financial support through a grant from the 
Simons Foundation (417314,MHW).

\bibliography{thesis_mlo}

\newcommand{\etalchar}[1]{$^{#1}$}
\begin{thebibliography}{BGK{\etalchar{+}}20}

\bibitem[AO20a]{AO20b}
Azam Asl and Michael~L. Overton.
\newblock {Analysis of limited-memory BFGS on a class of nonsmooth convex
  functions}.
\newblock {\em IMA Journal of Numerical Analysis}, 01 2020.
\newblock drz052.

\bibitem[AO20b]{AO20g}
Azam Asl and Michael~L. Overton.
\newblock {Analysis of the gradient method with an Armijo-Wolfe line search on
  a class of non-smooth convex functions}.
\newblock {\em Optimization Methods and Software}, 35(2):223--242, 2020.

\bibitem[Asl20]{AA20}
Azam Asl.
\newblock {\em Behavior of the Limited Memory BFGS Method on Nonsmooth
  Optimization Problems in Theory and Practice}.
\newblock PhD thesis, New York University, 2020.
\newblock
  https://cs.nyu.edu/media/publications/asl\_thesis\_final\_UtpoLsu.pdf.

\bibitem[BB88]{BB88}
Jonathan Barzilai and Jonathan~M. Borwein.
\newblock Two-point step size gradient methods.
\newblock {\em IMA Journal of Numerical Analysis}, 8(1):141--148, 1988.

\bibitem[BCL{\etalchar{+}}20]{BCLOS}
James~V. Burke, Frank~E. Curtis, Adrian~S. Lewis, Michael~L. Overton, and Lucas
  E.~A. Sim{\~o}es.
\newblock Gradient sampling methods for nonsmooth optimization.
\newblock In Adil~M. Bagirov, Manlio Gaudioso, Napsu Karmitsa, Marko~M.
  M{\"a}kel{\"a}, and Sona Taheri, editors, {\em Numerical Nonsmooth
  Optimization: State of the Art Algorithms}, pages 201--225. Springer
  International Publishing, Cham, 2020.

\bibitem[BGK{\etalchar{+}}20]{BGKMT20}
Adil~M. Bagirov, Manlio Gaudioso, Napsu Karmitsa, Marko~M. M{\"a}kel{\"a}, and
  Sona Taheri.
\newblock {\em Numerical Nonsmooth Optimization: State of the Art Algorithms}.
\newblock Springer International Publishing, Cham, 2020.

\bibitem[BKM14]{BKM14}
Adil Bagirov, Napsu Karmitsa, and Marko~M. M\"{a}kel\"{a}.
\newblock {\em Introduction to Nonsmooth Optimization}.
\newblock Springer, Cham, 2014.
\newblock Theory, Practice and Software.

\bibitem[BLO05]{BLO}
James~V. Burke, Adrian~S. Lewis, and Michael~L. Overton.
\newblock A robust gradient sampling algorithm for nonsmooth, nonconvex
  optimization.
\newblock {\em SIAM J. Optim.}, 15(3):751--779, 2005.

\bibitem[BV04]{BV}
Stephen Boyd and Lieven Vandenberghe.
\newblock {\em Convex Optimization}.
\newblock {Cambridge University Press}, March 2004.

\bibitem[CL20]{CurLi20}
Frank~E. Curtis and Minhan Li.
\newblock Gradient sampling methods with inexact subproblem solutions and
  gradient aggregation, 2020.
\newblock arXiv:2005.07822.

\bibitem[CMO17]{CMO17}
Frank~E. Curtis, Tim Mitchell, and Michael~L. Overton.
\newblock A {BFGS-SQP} method for nonsmooth, nonconvex, constrained
  optimization and its evaluation using relative minimization profiles.
\newblock {\em Optimization Methods and Software}, 32(1):148--181, 2017.

\bibitem[CO12]{CurOve12}
Frank~E. Curtis and Michael~L. Overton.
\newblock A sequential quadratic programming algorithm for nonconvex, nonsmooth
  constrained optimization.
\newblock {\em SIAM J. Optim.}, 22(2):474--500, 2012.

\bibitem[d'A08]{AdA08}
Alexandre d'Aspremont.
\newblock Smooth optimization with approximate gradient.
\newblock {\em SIAM J. Optim.}, 19(3):1171--1183, 2008.

\bibitem[DYC{\etalchar{+}}19]{madel}
Lijun Ding, Alp Yurtsever, Volkan Cevher, Joel~A. Tropp, and Madeleine Udell.
\newblock {An optimal-storage approach to semidefinite programming using
  approximate complementarity}, February 2019.
\newblock arXiv:1902.03373.

\bibitem[GB08]{gb08}
Michael Grant and Stephen Boyd.
\newblock Graph implementations for nonsmooth convex programs.
\newblock In V.~Blondel, S.~Boyd, and H.~Kimura, editors, {\em Recent Advances
  in Learning and Control}, Lecture Notes in Control and Information Sciences,
  pages 95--110. Springer-Verlag Limited, 2008.
\newblock http://stanford.edu/~boyd/graph\_dcp.html.

\bibitem[GB14]{cvx}
Michael Grant and Stephen Boyd.
\newblock {CVX}: Matlab software for disciplined convex programming, version
  2.1.
\newblock http://cvxr.com/cvx, March 2014.

\bibitem[GL18]{GL18}
J.~Guo and A.~Lewis.
\newblock Nonsmooth variants of {P}owell's {BFGS} convergence theorem.
\newblock {\em SIAM Journal on Optimization}, 28(2):1301--1311, 2018.

\bibitem[Gse14]{Gset}
{The {U}niversity of {F}lorida sparse matrix collection: Gset group}.
\newblock http://www.cise.ufl.edu/research/sparse/matrices/Gset/index.html,
  March 2014.

\bibitem[GW95]{GW95}
Michel~X. Goemans and David~P. Williamson.
\newblock Improved approximation algorithms for maximum cut and satisfiability
  problems using semidefinite programming.
\newblock {\em J. ACM}, 42(6):1115--€"1145, November 1995.

\bibitem[HOR14]{helOv}
C.~Helmberg, M.L. Overton, and F.~Rendl.
\newblock The spectral bundle method with second-order information.
\newblock {\em Optimization Methods and Software}, 29(4):855--876, 2014.

\bibitem[HR00]{helRen}
C.~Helmberg and F.~Rendl.
\newblock A spectral bundle method for semidefinite programming.
\newblock {\em SIAM Journal on Optimization}, 10(3):673--696, 2000.

\bibitem[Kiw85]{KIW85}
Krzysztof~C. Kiwiel.
\newblock {\em Methods of descent for nondifferentiable optimization}, volume
  1133 of {\em Lecture Notes in Mathematics}.
\newblock Springer-Verlag, Berlin, 1985.

\bibitem[Lem75]{LEM75}
C.~Lemar{\'e}chal.
\newblock An extension of {D}avidon methods to non differentiable problems.
\newblock {\em Math. Programming Stud.}, (3):95--109, 1975.

\bibitem[LN89]{LN89}
Dong~C. Liu and Jorge Nocedal.
\newblock On the limited memory {BFGS} method for large scale optimization.
\newblock {\em Math. Programming}, 45(3, (Ser. B)):503--528, 1989.

\bibitem[LO13]{LO13}
Adrian~S. Lewis and Michael~L. Overton.
\newblock Nonsmooth optimization via quasi-{N}ewton methods.
\newblock {\em Math. Program.}, 141(1-2, Ser. A):135--163, 2013.

\bibitem[LZ15]{LZ15}
A.~S. Lewis and S.~Zhang.
\newblock Nonsmoothness and a variable metric method.
\newblock {\em J. Optim. Theory Appl.}, 165(1):151--171, 2015.

\bibitem[MRT12]{Mohri12}
Mehryar Mohri, Afshin Rostamizadeh, and Ameet Talwalkar.
\newblock {\em Foundations of Machine Learning}.
\newblock Adaptive Computation and Machine Learning. MIT Press, Cambridge, MA,
  2012.

\bibitem[NB01]{NB01}
Angelia Nedi\'{c} and Dimitri~P. Bertsekas.
\newblock Incremental subgradient methods for nondifferentiable optimization.
\newblock {\em SIAM J. Optim.}, 12(1):109--138, 2001.

\bibitem[Nes05]{Nes05}
Yu. Nesterov.
\newblock Smooth minimization of non-smooth functions.
\newblock {\em Math. Program.}, 103(1, Ser. A):127--152, 2005.

\bibitem[Nes16]{YN16}
Yu. Nesterov.
\newblock Private communication, 2016.
\newblock Les Houches, France.

\bibitem[NW06]{NW06}
J.~Nocedal and S.~J. Wright.
\newblock {\em Numerical Optimization}.
\newblock Springer, New York, 2nd edition, 2006.

\bibitem[Ove88]{MO88}
M.L. Overton.
\newblock On minimizing the maximum eigenvalue of a symmetric matrix.
\newblock {\em SIAM J. Matrix Anal. Appl.}, 9:256--268, 1988.

\bibitem[Pow76]{Pow76b}
M.~J.~D. Powell.
\newblock Some global convergence properties of a variable metric algorithm for
  minimization without exact line searches.
\newblock In {\em Nonlinear Programming}, pages 53--72, Providence, 1976. Amer.
  Math. Soc.
\newblock SIAM-AMS Proc., Vol. IX.

\bibitem[RFP10]{FazRech}
Benjamin Recht, Maryam Fazel, and Pablo~A. Parrilo.
\newblock Guaranteed minimum-rank solutions of linear matrix equations via
  nuclear norm minimization.
\newblock {\em SIAM Review}, 52(3):471--501, 2010.

\bibitem[Van19]{VandenbergheCourse}
Lieven Vandenberghe.
\newblock {Optimization methods for large-scale systems}, 2019.
\newblock http://www.seas.ucla.edu/\~{ }vandenbe/236C/lectures/smoothing.pdf,
  Lecture Notes for ECE236C.

\bibitem[XW17]{XW17}
Yuchen Xie and Andreas Waechter.
\newblock {On the convergence of BFGS on a class of piecewise linear non-smooth
  functions}, December 2017.
\newblock arXiv:1712.08571.

\end{thebibliography}
\bibliographystyle{alpha}

\end{document}